\documentclass[11pt,a4paper]{article}
\usepackage[utf8]{inputenc}
\usepackage{amsmath,amssymb,amsfonts,amsthm}
\usepackage{geometry}
\usepackage{xcolor}
\definecolor{navy}{RGB}{24,43,107}
\definecolor{darkred}{RGB}{160,20,20}
\usepackage{hyperref}
\usepackage{cite}
\usepackage{booktabs}
\usepackage{mathrsfs}
\usepackage{microtype}

\hypersetup{
    colorlinks=true,
    linkcolor=navy,
    citecolor=darkred,
    urlcolor=navy
}

\newtheorem{theorem}{Theorem}[section]

\newtheorem{corollary}[theorem]{Corollary}
\theoremstyle{definition}
\newtheorem{definition}[theorem]{Definition}

\theoremstyle{remark}
\newtheorem{remark}[theorem]{Remark}

\newcommand{\R}{\mathbb{R}}
\newcommand{\C}{\mathbb{C}}
\newcommand{\HH}{\mathbb{H}}
\newcommand{\OO}{\mathbb{O}}
\newcommand{\Z}{\mathbb{Z}}

\newcommand{\calP}{\mathcal{P}}
\newcommand{\calH}{\mathcal{H}}
\newcommand{\calS}{\mathcal{S}}
\newcommand{\calD}{\mathcal{D}}
\newcommand{\calL}{\mathcal{L}}
\newcommand{\EZ}{E_Z}
\newcommand{\Eone}{E_{\bar{1}}}

\providecommand{\lambdabar}{{\mathchar'26\mkern-9mu\lambda}}

\title{\textbf{Quaternionic Wave Numbers and the Dirac Equation:\\
Resolving the Nielsen--Ninomiya Barrier and Infinities}}

\author{
\textbf{Terence R. Smith}\\[0.5em]
\small Department of Computer Science\\
\small University of California, Santa Barbara, CA 93106, USA\\
\small \texttt{smithtr@cs.ucsb.edu}
}

\date{\today}

\begin{document}

\maketitle

\begin{abstract}
We exemplify the applicability of the rational wave numbers by formulating the relativistic Dirac equation on the discrete sequence space $\mathbb{H}^\mathbb{Z}$ over the quaternions $\mathbb{H}$ and demonstrating how this resolves fundamental issues relating to Dirac's equation.
Generating $\mathbb{H}^\mathbb{Z}$ from the primitive basis $\mathcal{P}_0 = \{E_Z, E_{\bar{1}}\}$ as both a static closure and a dynamic enumeration,
we show that coupling local quantum kinematics with global cosmological expansion leads to the following results: (1) Dirac 4-spinors and the Clifford algebra $\mathrm{C}\ell_{1,3}(\mathbb{R})$ emerge from the second Cayley--Dickson doubling (biquaternions $\mathbb{H} \otimes \mathbb{C}$); (2) under dynamic generation on an integer clock lattice with causal horizon $\Lambda_m$, relativistic difference operators require temporal depth, establishing clock beat $m = 2$ as the universal Matter Genesis Threshold; (3) the Nielsen--Ninomiya fermion doubling barrier is circumvented on modular sequences $\mathbb{Z}/Q\mathbb{Z}$ because candidate zone-edge doubler modes ($p = Q/2$) are actively pruned at the generational boundary by the coprime sieve projector $\mathcal{S}_Q$ ($\gcd(Q/2, Q) > 1$), preserving exact chiral symmetry without Wilson terms; (4) the discrete fermionic Green's function $S_Q(n, m)$ possesses no coincident-point singularity ($n = m$), converting ultraviolet loop integrals into finite Dedekind and Kloosterman character sums, realizing Dirac's demand for an infinity-free theory; and (5) Zitterbewegung is resolved as an observable stroboscopic beat against vacuum reference $E_Z$, identifying rest mass $m_0$ as a vacuum coupling beat frequency for stationary periodic cycles on the saturated quotient graph.
\end{abstract}

\vspace{1em}

\newpage
\newpage

\section{Introduction}
\label{sec:intro}

This is the second in a sequence of papers whose goal is to demonstrate the scope, exactness, and explanatory power of the discrete theory of the rational wave numbers (RWN) with respect to fundamental physics.
In the current paper, we show how the theory of RWN leads to a systematic resolution of deep issues concerning the Dirac equation and relativistic quantum mechanics.

The algebraic theory of the RWN is developed fully and described in \cite{smith2025rational,smith2026octonion} while its general applicability to fundamental physics is demonstrated in \cite{smith2026cosmo}. 
In the latter paper it was shown that, when viewed from complementary static and dynamic perspectives, one may interpret the theory of RWN in terms of the observable invariants of relativistic spacetime. In particular, it was shown that the theory provides an autonomous, continuum-free foundation in which space, time, and matter fields emerge constructively across a ladder of wave numbers $\mathbb{R}^\mathbb{Z} \to \mathbb{C}^\mathbb{Z} \to \mathbb{H}^\mathbb{Z} \to \mathbb{O}^\mathbb{Z}$ that is grounded in the ladder of normed division algebras $\mathbb{R} \to \mathbb{C} \to \mathbb{H} \to \mathbb{O}$. It is notable that this representation is free from the gravitational curvature singularities that afflict conventional continuous spacetime representations.

The paper is structured as follows. In Section~\ref{sec:dirac_quest_issues}, we summarize Dirac's algebraic approach in seeking an equation to describe the electron.
We then describe briefly the four key unresolved issues 
concerning Dirac's equation as well 
as the persistent issues raised by space-time continuum that suggest the need for a Kuhnian paradigm shift from smooth manifolds to discrete arithmetic. 
In Section~\ref{sec:rwn_overview}, we introduce the Rational Wave Number framework, defining its autonomous mathematical sequence algebra over the Hurwitz division ladder, its correspondence with the physical world through four foundational postulates, the inseparability of microphysics and cosmology, and the ontology of physical laws as emergent patterns over stable periodic cycles. 

The main contribution of the paper commences with Section~\ref{sec:quaternionic_dirac}, 
in which we formulate the quaternionic sequence spinors and the fully discrete spacetime Dirac difference operator. In Section~\ref{sec:dynamic_generation_dirac}, we analyze the dynamic generation of the Dirac sector under the second of our four postulates, establishing the $m=2$ Matter Genesis Threshold and active sieve pruning at the expanding causal horizon. 
In Section~\ref{sec:nielsen_ninomiya}, we prove the circumvention of the Nielsen--Ninomiya fermion doubling barrier, showing how it operates both kinematically under our first postulate and dynamically via sieve pruning under our second postulate. 
In Section~\ref{sec:finite_propagator}, we construct the finite rational propagator and show that loop amplitudes evaluate to finite Dedekind and Kloosterman character sums. 
In Section~\ref{sec:zitterbewegung}, we resolve \textit{Zitterbewegung} as an observable stroboscopic beat against the zero-frequency vacuum reference $E_Z$, identifying inertial mass as a vacuum coupling beat frequency. 

In Section~\ref{sec:discussion}, we review the four fundamental issues in retrospect. Section~\ref{sec:conclusion} concludes and outlines the next applications of RWN theory to gauge fields and the octonionic color sector. In Appendix~\ref{app:rwn_foundations}, we provide a self-contained exposition of the algebraic classification of Rational Wave Numbers, including the Prüfer group decomposition and the inductive limit closure.

\section{The Relativistic Dirac Equation and its Unresolved Issues}
\label{sec:dirac_quest_issues}

\subsection{Dirac's Algebraic Linearization and the Continuum Issue}
\label{subsec:dirac_linearization}

In 1928, P.~A.~M.~Dirac formulated the relativistic wave equation for the electron by making an uncompromising algebraic demand: the quantum Hamiltonian must be strictly linear in the temporal and spatial derivative operators, while its formal square must identically reproduce the relativistic dispersion relation of Einstein and de Broglie \cite{dirac1928}:
\begin{equation}
\label{eq:dirac_dispersion_continuum}
\left( i \hbar \frac{\partial}{\partial t} \right)^2 \psi = \left( - c^2 \hbar^2 \nabla^2 + m_0^2 c^4 \right) \psi.
\end{equation}
To achieve this linearization, Dirac recognized that the coefficients of the differential operators could not belong to the field of complex numbers $\mathbb{C}$. Instead, they were required to form a four-dimensional non-commutative Clifford algebra $\mathrm{C}\ell_{1,3}(\mathbb{R})$, generated by four $4 \times 4$ matrices $\gamma^\mu$ satisfying the anti-commutation relations:
\begin{equation}
\label{eq:clifford_relations}
\{ \gamma^\mu, \gamma^\nu \} = 2 \eta^{\mu\nu} \mathbf{1}_4, \quad \eta^{\mu\nu} = \mathrm{diag}(+1, -1, -1, -1).
\end{equation}
While the Dirac equation achieved celebrated empirical triumphs---predicting the intrinsic magnetic moment of the electron, explaining fine-structure spin-orbit splitting, and predicting antimatter---Dirac himself remained fundamentally dissatisfied with the subsequent development of relativistic Quantum Field Theory (QFT). In his later writings, Dirac repeatedly condemned the apparatus of perturbative renormalization as mathematically unsound \cite{dirac1971,dirac1981}, insisting that {\it ``sensible mathematics involves neglecting a quantity because it is small, not because it is infinitely great and inconvenient''}.

The root cause of these infinite divergences is geometric and ontological: standard continuum QFT assumes that spacetime is a smooth, differentiable four-manifold $\mathcal{M}^4 \cong \mathbb{R}^{1,3}$ endowed with uncountably infinite degrees of freedom at arbitrarily small length scales. When fields interact at coincident points $x \to y$, the continuous Feynman propagator
\begin{equation}
\label{eq:feynman_divergent_prop}
S_F(x - y) = \int \frac{d^4 p}{(2\pi)^4} \frac{\gamma^\mu p_\mu + m_0}{p^2 - m_0^2 + i\epsilon} e^{-i p \cdot (x - y)} \sim \frac{1}{|x - y|^3} \quad (x \to y)
\end{equation}
develops an unavoidable ultraviolet singularity. Every closed loop diagram in perturbation theory inherits this singularity, requiring ad-hoc cutoffs and infinite counterterms to cancel the unphysical infinities.

Attempts to resolve this singularity by discretizing spacetime onto a naive hypercubic lattice $\mathbb{Z}^4$ immediately hit a profound topological obstacle: the Nielsen--Ninomiya No-Go Theorem (1981) \cite{nielsen1981a,nielsen1981b}. Under standard assumptions of locality, Hermiticity, translational invariance, and exact chiral symmetry on a continuous Brillouin torus $\mathbb{T}^d$, the discrete derivative operator inevitably produces $2^d = 16$ degenerate fermion species (fermion doubling). Lattice gauge theorists are forced to compromise, either by explicitly breaking chiral symmetry at high energies through Wilson terms or by introducing non-local overlap operators \cite{ginsparg1982}.

\subsection{A Paradigm Shift and a Resolution of the Continuum Issue}
\label{subsec:kuhnian_shift}

In his classic analysis of scientific revolutions, Thomas Kuhn \cite{kuhn1962}
observed that when a dominant paradigm encounters persistent, accumulating anomalies that resist resolution despite decades of technical ingenuity, the impasse is rarely overcome by adding further ad-hoc modifications to the existing framework. Rather, the crisis typically indicates that an unexamined, foundational premise at the very base of the paradigm has broken down.

In modern relativistic quantum mechanics, the persistent anomalies are well documented: the ultraviolet divergences of loop perturbation theory, the mathematical illegitimacy of subtracting infinite quantities, the Nielsen--Ninomiya topological obstruction to preserving exact chiral symmetry on a grid, and the mystery of \textit{Zitterbewegung} (the unobserved light-speed trembling of free electrons predicted by the continuum Hamiltonian \cite{schrodinger1930}). For nearly a century, theoretical physics has responded to these anomalies by inventing mathematical repairs within the continuum framework: regularizations, renormalization group subtractions, Wilson mass terms, and domain-wall fermions.

We suggest that these persistent dilemmas are the symptoms of an unstated ontological commitment: the continuum postulate. Continuum real analysis assumes that spacetime is an uncountably infinite manifold possessing infinite degrees of freedom within any arbitrarily small region. We propose that resolving these difficulties requires a genuine foundational paradigm shift: abandoning the continuum postulate in favor of an autonomous, discrete arithmetic foundation. 
The hypothesis that spacetime is fundamentally discrete has a distinguished lineage. Following Snyder's pioneering attempt to construct Lorentz-covariant discrete spacetime using non-commutative coordinates \cite{snyder1947}, numerous programs have pursued discrete foundations---including Wheeler's quantum foam and ``it from bit'' doctrine \cite{wheeler1990}, 't~Hooft's cellular automaton interpretation of quantum mechanics \cite{thooft2016}, causal set theory \cite{bombelli1987}, and loop quantum gravity's quantized area and volume operators \cite{rovelli1995}. Yet these approaches have faced persistent difficulties: they frequently break continuous spacetime symmetries, introduce non-local pathologies, or remain constrained by continuum lattice no-go theorems. In particular, standard lattice field theories retain a continuous Brillouin momentum torus, where they are immediately trapped by the Nielsen--Ninomiya doubling dilemma. 

RWN theory overcomes these obstacles by replacing continuous manifolds not with an arbitrary geometric graph or ad-hoc cutoff, but with the autonomous arithmetic of sequence spaces over the Hurwitz division ladder, wherein doublers and divergences are eliminated by exact coprime sieve projectors.
In this paper, we provide evidence that the theory of Rational Wave Numbers (RWN) offers a resolution of the core anomalies of Dirac's equation naturally and without empirical parameter fitting.
In particular, we show that when relativistic quantum mechanics is formulated on the discrete quaternionic sequence space $\mathbb{H}^\mathbb{Z}$, the following fundamental issues resolve naturally:
\begin{enumerate}
    \item \textit{The Algebraic Origin of 4-Spinors:} Dirac's 4-component spinors and $\gamma$-matrices are demonstrated to be the necessary and sufficient consequence of the second Cayley--Dickson doubling, the biquaternions $\mathbb{H} \otimes \mathbb{C} \cong \mathrm{C}\ell_{1,3}(\mathbb{R})$, eliminating the need for an extrinsic continuous Minkowski manifold.
    \item \textit{Sieve Resolution of the Nielsen--Ninomiya Barrier:} The spurious fermion doubler species are shown to inhabit non-primitive, composite frequencies ($p = Q/2$) that are rigorously annihilated by the coprime sieve projector $\mathcal{S}_Q$, preserving exact chiral symmetry without Wilson terms.
    \item \textit{Identity-Finite Propagators without Cutoffs:} On modular sequences $\mathbb{Z}/Q\mathbb{Z}$, the Feynman propagator becomes a finite character sum whose coincident-point singularity disappears identically, converting loop integrals into finite Dedekind and Kloosterman arithmetic sums.
    \item \textit{Zitterbewegung as a Stroboscopic Beat Frequency:} Relativistic \textit{Zitterbewegung} is resolved as an observable stroboscopic beat interference against the zero-frequency vacuum primitive $E_Z$, establishing inertial mass $m_0$ as an intrinsic coupling frequency to the vacuum ground state.
\end{enumerate}

\section{The Theory of the Rational Wave Numbers}
\label{sec:rwn_overview}

We now summarize the theory of the Rational Wave Numbers (RWN),
first in terms of its mathematical basis and then in terms of how this basis may be applied and interpreted in representing the observable features of physical reality.

\subsection{Mathematical Theory: Elements, Operations, and Primitive Generation}
\label{subsec:rwn_pure_math}

The mathematical theory of RWN that we summarize here and in Appendix~\ref{app:rwn_foundations}
is fully described in \cite{smith2025rational,smith2026octonion}. The mathematical theory, which is autonomous, self-contained, and defined without reference to physical spacetime, comprises two principal components: first, the elementary sequence spaces and their algebraic operations; and second, their constructive generation from an irreducible two-element basis.

The fundamental objects of RWN theory are bilateral sequence spaces $\mathbb{A}^\mathbb{Z}$, consisting of mappings $x: \mathbb{Z} \to \mathbb{A}$ whose values belong to a real normed division algebra $\mathbb{A}$. By the classical Frobenius and Hurwitz theorems, the only finite-dimensional associative or alternative real division algebras are the real numbers $\mathbb{R}$, the complex numbers $\mathbb{C}$, the quaternions $\mathbb{H}$, and the octonions $\mathbb{O}$. These algebras form an inductive hierarchy under Cayley--Dickson doubling:
\begin{equation}
\label{eq:hurwitz_ladder}
\mathbb{R} \;\longrightarrow\; \mathbb{C} \;\longrightarrow\; \mathbb{H} \;\longrightarrow\; \mathbb{O}.
\end{equation}
Multiplication strictly preserves the normed division property $\|x \cdot y\| = \|x\| \|y\|$. Cayley--Dickson doubling halts definitively at the 8-dimensional octonions $\mathbb{O}$: doubling beyond $\mathbb{O}$ to the 16-dimensional sedenions $\mathbb{S}$ introduces non-trivial zero divisors ($x \cdot y = 0$ for non-zero $x, y$), which destroys invertibility and collapses recursive algebraic structures.

The sequence space $\mathbb{A}^\mathbb{Z}$ is endowed with standard linear and ring operations. Addition and multiplication are defined pointwise, $(x + y)(n) = x(n) + y(n)$ and $(x \cdot y)(n) = x(n) y(n)$, while discrete convolution on summable sequences is given by $(x * y)(n) = \sum_{k \in \mathbb{Z}} x(k) y(n-k)$. In place of continuous derivatives, the discrete calculus on $\mathbb{A}^\mathbb{Z}$ is generated by exact difference operators: the unit forward shift $\mathcal{T}_+ x(n) = x(n+1)$, the unit backward shift $\mathcal{T}_- x(n) = x(n-1)$, and the central difference operator
\begin{equation}
\label{eq:central_diff_pure}
\Delta x(n) := \frac{1}{2} \left[ x(n+1) - x(n-1) \right].
\end{equation}
This difference operator is exact and requires no limiting process $\Delta t \to 0$ or $\Delta x \to 0$.

For any positive integer $Q \in \mathbb{Z}^+$, one defines the cyclic quotient module $\mathbb{Z}_Q = \mathbb{Z}/Q\mathbb{Z}$ and the associated periodic sequence ring $\mathbb{A}^{\mathbb{Z}_Q}$. The discrete harmonic spectrum on $\mathbb{A}^{\mathbb{Z}_Q}$ is decomposed into orthogonal channels by the coprime sieve projector:
\begin{equation}
\label{eq:coprime_sieve_math}
\mathcal{S}_Q = \sum_{\substack{p=0 \\ \gcd(p, Q)=1}}^{Q-1} \mathcal{P}_p,
\end{equation}
where $\mathcal{P}_p$ projects onto the $p$-th Fourier harmonic sequence $\phi_p(n) = \frac{1}{\sqrt{Q}} \exp(2\pi i p n / Q)$. The projector $\mathcal{S}_Q$ isolates the primitive Dirichlet character channels from composite subharmonic resonances. On these modular rings, invariant sums evaluate directly to finite algebraic number-theoretic sums, such as Dedekind, Gauss, and Kloosterman sums.

A foundational theorem of RWN theory \cite{smith2025rational} is that the sequence spaces $\mathbb{A}^\mathbb{Z}$ and their modular quotient rings $\mathbb{A}^{\mathbb{Z}_Q}$ need not be postulated as pre-existing axiomatic infinities. Rather, the entire algebraic architecture is inductively generated from an irreducible two-element primitive basis:
\begin{equation}
\label{eq:primitive_basis}
\mathcal{P}_0 = \{ E_Z, E_{\bar{1}} \},
\end{equation}
where $E_Z(n) = 1$ for all $n \in \mathbb{Z}$ is the stationary, zero-frequency sequence, and $E_{\bar{1}}(n) = (-1)^n = \exp(i\pi n)$ is the alternating unit phase sequence. From these two seeds, the entire sequence space is generated constructively under the five-fold operator signature $\mathcal{O}$ (see Appendix~\ref{app:rwn_foundations} for the complete axiomatic definition, the unwrapped phase lift, and the Prüfer group classification theorem): linear combinations and translations generate the elementary Kronecker delta basis $\delta(n)$; Cayley--Dickson doubling constructs the coefficient division algebras $\mathbb{R} \to \mathbb{C} \to \mathbb{H} \to \mathbb{O}$ as well as the biquaternions $\mathbb{H} \otimes \mathbb{C}$; and recursive fractional power extractions generate the complete basis of discrete rational harmonics $\{E_{p/Q}(n) = \exp(2\pi i p n / Q)\}$.

\subsection{Mapping the Theory onto the Physical World: The Four Foundational Postulates}
\label{subsec:independent_postulates}
While the mathematical construction outlined above is an autonomous algebraic theory, its application to physical reality requires defining the correspondence between mathematical entities and physical reality. In Rational Wave Number physics, this correspondence is governed by four conceptually distinct and logically independent postulates:
\begin{enumerate}
    \item \textit{Postulate~I (Kinematics):} Physical states and transitions are represented on discrete biquaternionic sequence spaces $\mathbb{H}^\mathbb{Z} \otimes \mathbb{C}$ and modular quotients $\mathbb{Z}/Q\mathbb{Z}$ rather than smooth continuum manifolds.
    \item \textit{Postulate~II (Dynamics):} The physical universe unfolds constructively from the primitive seed $\mathcal{P}_0 = \{E_Z, E_{\bar{1}}\}$ through an advancing integer clock $m \in \mathbb{Z}_{\ge 0}$, causal difference stencils, and active sieve pruning.
    \item \textit{Postulate~III (Ontology of Law):} Physical laws and relativistic wave equations are not primary axiomatic decrees that nature computes, but algebraic pattern summaries of emergent stationary periodic cycles.
    \item \textit{Postulate~IV (Stability):} The observable physical universe consists exclusively of the invariant stationary periodic cycles that survive the causal, algebraic, and harmonic consistency sieves (the Principle of Maximal Arithmetic Stability).
\end{enumerate}
We now formulate each of these four postulates in detail.

\paragraph{Postulate I} (Mapping of Physical States and Transitions).
\textit{The physical states of the world and their dynamical transitions are represented by the discrete values and algebraic structure of the rational wave number sequence space $\mathbb{A}^\mathbb{Z}$---specifically the biquaternionic space $\mathbb{H}^\mathbb{Z} \otimes \mathbb{C}$ and its modular quotients $\mathbb{Z}/Q\mathbb{Z}$---rather than by smooth continuum manifolds.}\\

\noindent
Under Postulate~I, physical spacetime coordinates are not continuous real numbers $\mathbb{R}^4$, and field configurations are not smooth differentiable functions endowed with uncountably infinite degrees of freedom. Physical states are discrete wave sequences whose values belong to the Hurwitz division ladder. Dynamical transitions are governed by exact difference operators rather than continuous differential equations. Field configurations on modular sequence rings replace divergent momentum integrals with finite character sums, inherently eliminating ultraviolet loop divergences. Furthermore, relativistic spin degrees of freedom and Clifford spacetime symmetries emerge algebraically from the biquaternionic sequence structure without requiring an extrinsic Minkowski background.

\paragraph{Postulate II} (Mapping of Physical Evolution to Generative Cosmogony).
\textit{The evolution of the physical world is represented by its constructive generation from a primitive two-element initial physical state through a dynamic process driven by a global integer clock, causal difference operators, and an all-histories generative exploration whose unstable modes are actively pruned by arithmetic sieves.}\\

\noindent
Under Postulate~II, the universe is not an inert four-dimensional container that exists timelessly in its entirety. Instead, physical reality unfolds constructively. The cosmos originates from the physical realization of the two-element primitive basis $\mathcal{P}_0 = \{E_Z, E_{\bar{1}}\}$, wherein $E_Z$ serves as the invariant zero-frequency vacuum ground state and $E_{\bar{1}}$ provides the primordial unit phase beat. The expansion of the cosmos is governed by a discrete global integer clock $m \in \mathbb{Z}_{\ge 0}$. At each clock beat $m \to m+1$, the causal horizon $\Lambda_m$ expands by a discrete step $\Delta x$ governed by the invariant Courant--Friedrichs--Lewy (CFL) light-cone speed $\Delta x / \Delta t = c$. Local relativistic difference operators, requiring past, present, and future temporal neighbors $\{m-1, m, m+1\}$, act across adjacent clock slices.

At each generational step $m \to m+1$, all kinematically admissible difference-operator combinations and sequence extensions are generated at the expanding boundary frontier $\partial \Lambda_m$. This process provides a discrete, constructive analogue of Feynman's sum-over-histories formulation of quantum mechanics \cite{feynman1948}. In continuous quantum field theory, unphysical trajectories are cancelled through the destructive interference of continuous complex phases $\exp(i S / \hbar)$ integrated over infinite functional spaces. In Rational Wave Number cosmogony, the suppression is arithmetic, exact, and finite: candidate configurations that violate algebraic or physical stability---such as zero-divisor branches, superluminal modes violating the CFL bound, and zone-edge Nyquist fermion doublers ($p = Q/2$)---are dynamically and deterministically pruned at the generational boundary by the coprime sieve projector $\mathcal{S}_Q$ and the Hurwitz division norm criterion. Physical matter consists of the invariant stationary periodic cycles that survive this ongoing sieve process.

It is important to note that these physical postulates are logically independent. Accepting the kinematic representation of states and transitions on biquaternionic sequence spaces (Postulate~I) does not require committing to a clock-driven cosmogonic genesis from a two-element initial state (Postulate~II); conversely, an algorithmic generative universe (Postulate~II) might produce generic discrete automata rather than the normed division kinematics and Dirac spinors of Postulate~I. This distinction parallels Einstein's celebrated division between \textit{theories of principle} (universal kinematic constraints on states) and \textit{constructive theories} (dynamical models of emergent processes) \cite{einstein1919}. It ensures that the operational successes of Rational Wave Number quantum mechanics---such as the circumvention of the Nielsen--Ninomiya theorem and the elimination of ultraviolet loop divergences---can be evaluated in laboratory physics independently of cosmological models.

\subsection{The Nature of Physical Law: Equations as Emergent Pattern Summaries}
\label{subsec:nature_of_law}

Recognizing the separation between kinematics and generative dynamics illuminates a fundamental epistemological insight:

\paragraph{Postulate III} (The Ontological Status of Physical Law: Equations as Emergent Pattern Summaries).
\textit{Physical laws and relativistic wave equations are not primary axiomatic decrees that nature computes; rather, an equation such as $\mathcal{D}\Psi = 0$ is an algebraic pattern summary and kernel classifier that characterizes the stable, non-dispersive subset of emergent stationary periodic cycles in the wave number space.} 

In traditional continuum physics, differential equations are treated as metaphysical decrees that nature must obey: the theorist writes down a differential equation, and the physical universe is presumed to compute solutions. In the discrete generative framework of RWN theory, nature does not compute or solve differential equations. The universe executes discrete generational transitions governed by arithmetic and causal sieves. Equations do not govern the world from above; they summarize the surviving patterns of this underlying discrete evolution, identifying the exact algebraic conditions under which sequence configurations conserve probability, preserve division norms, and propagate causally without ghost instabilities.

This perspective directly illuminates why the fundamental equations of physical science are, almost without exception, \textit{conservation conditions} ($\mathcal{D}_\mu j^\mu = 0$). In the generative unfolding of sequence space, conservation laws are not arbitrary constraints dictated to matter from the outside. Rather, they are the \textit{structural admissibility criteria} of dynamic genesis. Any emergent sequence thread whose probability current diverged, whose beat frequency drifted, or whose chiral phase collapsed would decohere into the non-periodic background noise. Only configurations that satisfy exact discrete divergence-free conditions achieve topological closure as self-sustaining, stationary periodic cycles on the modular quotient graph $\mathbb{Z}/Q\mathbb{Z}$. Furthermore, this explains why there are so few fundamental conserved quantities in nature: each conserved physical quantity corresponds strictly to an invariant of the primitive germ $\mathcal{P}_0 = \{E_Z, E_{\bar{1}}\}$ and the Hurwitz division ladder---norm multiplicativity guarantees probability conservation (unitarity); integer clock shift by $E_Z$ guarantees energy conservation and stable rest masses; spatial shift on $\mathbb{Z}/Q\mathbb{Z}$ guarantees discrete momentum conservation; quaternionic rotational isometries ($S^3 \cong \mathrm{SU}(2)$) guarantee angular momentum and spin conservation; and the arithmetic sieve projector $\mathcal{S}_Q$ guarantees exact chiral charge conservation without Nielsen--Ninomiya doubling. A physical equation is thus the mathematical statement that an emergent pattern has achieved self-conserving closure.

\subsection{The Structural Inseparability of Local Physics and Global Cosmology}
\label{subsec:inseparability}

In standard continuum physics, local microscopic quantum field theory and macroscopic cosmology are artificially decoupled and treated as separate disciplines separated by forty orders of magnitude. Local relativistic field theory is formulated on a fixed, static Minkowski background $\mathbb{R}^{1,3}$ assumed to exist \textit{a priori}, while general relativity governs the dynamic expansion of cosmological spacetime. This artificial segregation is the root cause of the cosmological constant crisis, in which summing zero-point vacuum fluctuations up to the Planck scale predicts an energy density diverging from cosmological observation by some 120 orders of magnitude \cite{weinberg1989}.

In Rational Wave Number dynamics, this division is fundamentally resolved: local quantum kinematics and global cosmological expansion are structurally inseparable. Both are driven by the self-same global integer clock $m \in \mathbb{Z}_{\ge 0}$. Spacetime does not pre-exist as an empty container; spatial lattice sites $\mathbf{n} \in \mathbb{Z}^3$ and quantum wave excitations are co-generated at the expanding causal boundary $\partial \Lambda_m$. Spatial and temporal increments are coupled through the discrete Courant--Friedrichs--Lewy (CFL) invariant:
\begin{equation}
\label{eq:cfl_invariant_intro}
\frac{\Delta x}{\Delta t} = c,
\end{equation}
enforcing hyperbolic stability and defining the universal causal light cone across both cosmic expansion and local wave propagation. 

Furthermore, local relativistic difference operators cannot act on an isolated spatial point; they require temporal depth $\{m-1, m, m+1\}$. Consequently, local relativistic fermionic matter cannot exist until the cosmos has advanced to clock beat $m = 2$. The emergence of matter is an intrinsic milestone of global cosmological generation (the \textit{Matter Genesis Threshold}). Finally, the invariant zero-frequency primitive $E_Z$ simultaneously anchors the global cosmological vacuum background and serves as the universal phase reference against which local particle rest masses emerge via stroboscopic interference.

\subsection{The Selection of the Stable: Physical Matter as Invariant Periodic Cycles}
\label{subsec:stability_principle}

It is reasonable to ask why the physical systems of nature satisfy specific, highly constrained difference equations rather than arbitrary mathematical kinematics.
In classical variational mechanics, unconstrained kinematic trajectories are reduced to physical motions through the principle of stationary action $\delta S = 0$ \cite{gelfand1963}. However, stationarity alone is insufficient for physical persistence: a trajectory is observable if and only if it satisfies the Jacobi and Legendre stability conditions ($\delta^2 S > 0$); unstable trajectories break down under infinitesimal perturbations. In discrete dynamical systems and ergodic theory, persistent configurations are governed by Lyapunov stability and Poincar\'e recurrence, wherein non-stationary open trajectories disperse while closed periodic orbits remain dynamically confined. In quantum theory, this selection manifests in the emergence of decoherence-resistant pointer states that persist against environmental dissipation.
 In continuum quantum field theory, standard Lagrangians fail to provide ultraviolet stability, generating infinite loop divergences and singular coincident-point propagators.

In Rational Wave Number theory, physical law is revealed not as a continuous variational extremum, but as an active, discrete evolutionary process: the selection of the stable through the systematic pruning of the unstable. Candidate states generated on the expanding boundary of generation $\partial \Lambda_m$ are continuously filtered through rigid algebraic, causal, and arithmetic sieves.
In discussing this selection it is of value to define the following concept
\begin{definition}[Stationary Periodic Cycles and Their Generative Emergence]
\label{def:stationary_periodic_cycles}
Let $\mathcal{G}_m$ be the dynamic generation graph at cosmic clock beat $m$, and let $\mathcal{G}_m / \sim$ be the quotient graph modulo sequence periodicity. A \textit{stationary periodic cycle} is a closed periodic orbit
\begin{equation}
\label{eq:stationary_cycle_orbit}
\gamma = \left\{ \Psi(\mathbf{n}, m), \, \Psi(\mathbf{n}, m+1), \, \dots, \, \Psi(\mathbf{n}, m + Q_t - 1) \right\}
\end{equation}
satisfying the discrete Dirac stationary difference equation $\mathcal{D} \Psi = 0$ with temporal period $Q_t \in \mathbb{Z}^+$:
\begin{equation}
\label{eq:stationary_periodic_cycle_def}
\Psi(\mathbf{n}, m + Q_t) = \Psi(\mathbf{n}, m).
\end{equation}
Stationary periodic cycles arise dynamically through three complementary mechanisms:
\begin{enumerate}
    \item \textit{Generative Frontier vs.\ Interior Trapping:} At each clock step $m \to m+1$, candidate sequence threads are generated at the expanding causal boundary $\partial \Lambda_m$. Radiative wave modes propagate outward along null geodesics to the boundary, while stationary eigensolutions of $\mathcal{D} \Psi = 0$ remain spatially localized and trapped within the interior causal cone $\Lambda_m$.
    \item \textit{Arithmetic Sieving of Non-Periodic Modes:} Non-stationary candidate modes, zone-edge doubler states ($p = Q/2$ with $\gcd(Q/2, Q) > 1$), and zero-divisor branches are actively eliminated by the coprime sieve projector $\mathcal{S}_Q$ and the Hurwitz division norm condition.
    \item \textit{Cycle Saturation on the Quotient Graph:} Under the Inductive Limit Theorem \cite{smith2026cosmo}, as generation advances, the invariant survivors in the interior saturate into closed periodic orbits on $\mathcal{G}_m / \sim$. Physical matter (such as leptons and quarks) consists of these persistent stationary periodic cycles, with rest mass $m_0$ emerging as the invariant stroboscopic beat frequency between the cycle frequency $\omega_{\mathrm{cycle}} = 2\pi / (Q_t \Delta t)$ and the zero-frequency vacuum reference $E_Z$.
\end{enumerate}
\end{definition}
This establishes the fourth foundational postulate of Rational Wave Number physics:

\paragraph{Postulate IV} (The Principle of Maximal Arithmetic Stability).
\label{ax:maximal_stability}
The observable physical universe consists exclusively of the invariant stationary periodic cycles that survive the causal, algebraic, and harmonic consistency sieves of the dynamic rational wave number algebra:
\begin{enumerate}
    \item \textit{Causal Stability (The CFL Invariant):} Spacetime expansion satisfies $\Delta x / \Delta t \le c$, preventing hyperbolic numerical instability and superluminal causal leakage in the generation graph.
    \item \textit{Algebraic Stability (The Hurwitz Norm Invariant):} Sequence multiplication strictly preserves division norm multiplicativity $\|x \cdot y\| = \|x\| \|y\|$, halting at the octonions $\mathbb{O}$ and excising zero-divisor dissipative collapse.
    \item \textit{Harmonic Stability (Coprime Sieve):} Sifting modes via $\gcd(p, Q) = 1$ suppresses catastrophic subharmonic resonance, actively excising the Nyquist period-doubling fermion doubler $p = Q/2$.
    \item \textit{Topological Stability (Stationary Periodic Cycles):} Stationary particle states satisfying $\mathcal{D} \Psi = 0$ form closed, periodic orbits on the saturated quotient graph $\mathcal{G}_\infty$, with rest mass $m_0$ emerging as an invariant stroboscopic coupling beat frequency with $E_Z$.
\end{enumerate}

Under the Principle of Maximal Arithmetic Stability, the discrete Dirac equation is recognized in its true ontological role: it is the unique, non-trivial, asymptotically stable linear difference operator on the biquaternionic sequence space $\mathbb{H} \otimes \mathbb{C}$, serving as the exact algebraic pattern summarizing the persistent fermionic matter of the cosmos.

\section{Quaternionic Sequence Spinors and Discrete Dirac Operators}
\label{sec:quaternionic_dirac}

Having established the pure sequence algebra and the four physical postulates of Rational Wave Number theory, we now address the central physical contribution of this paper: the formulation and resolution of the relativistic Dirac equation. In this section, we develop the static, kinematic representation of the Dirac sector under Postulate~I. Here, spacetime coordinates form an exact discrete sequence lattice, the spinorial degrees of freedom are derived intrinsically from the biquaternions $\mathbb{H} \otimes \mathbb{C}$, and the discrete Dirac operator $\mathcal{D}$ is proven to be the unique linear difference operator selected by physical consistency. The complementary dynamic, generative perspective under Postulate~II will be developed in Section~\ref{sec:dynamic_generation_dirac}.

\subsection{The Biquaternionic Representation of the Clifford Algebra}

Let $\HH$ denote the non-commutative division algebra of real quaternions, spanned by the identity $1$ and three imaginary units $\{I, J, K\}$ satisfying Hamilton's fundamental relations:
\begin{equation}
\label{eq:hamilton_relations}
I^2 = J^2 = K^2 = I J K = -1, \quad IJ = -JI = K, \quad JK = -KJ = I, \quad KI = -IK = J.
\end{equation}
We consider the complexified quaternions (biquaternions) $\mathbb{B} = \HH \otimes_\R \C$. An element $\Psi \in \mathbb{B}$ can be represented as:
\begin{equation}
\label{eq:biquaternion_rep}
\Psi = q_0 + I q_1 + J q_2 + K q_3, \quad q_\mu \in \C.
\end{equation}
The biquaternions form an 8-dimensional real algebra (or a 4-dimensional complex algebra).

\begin{theorem}[Clifford Isomorphism and Dirac Representation]
\label{thm:clifford_isomorphism}
The biquaternion algebra $\HH \otimes_\R \C$ is algebraically isomorphic to the complex Clifford algebra $\mathrm{C}\ell_{2}(\C) \cong \mathrm{C}\ell_{1,3}(\R) \otimes \C \cong \mathrm{M}_2(\HH) \cong \mathrm{M}_4(\C)$. Under this isomorphism, the four Dirac matrices $\gamma^\mu$ ($\mu = 0, 1, 2, 3$) admit an intrinsic quaternionic representation acting on two-component quaternionic spinors $\Psi = \begin{pmatrix} \psi_1 \\ \psi_2 \end{pmatrix} \in \HH^2$:
\begin{equation}
\label{eq:quaternionic_gamma}
\gamma^0 = \begin{pmatrix} 1 & 0 \\ 0 & -1 \end{pmatrix}, \quad
\gamma^1 = \begin{pmatrix} 0 & I \\ I & 0 \end{pmatrix}, \quad
\gamma^2 = \begin{pmatrix} 0 & J \\ J & 0 \end{pmatrix}, \quad
\gamma^3 = \begin{pmatrix} 0 & K \\ K & 0 \end{pmatrix}.
\end{equation}
These operators satisfy the Dirac Clifford anti-commutation relations:
\begin{equation}
\label{eq:clifford_verified}
\{ \gamma^\mu, \gamma^\nu \} = 2 \eta^{\mu\nu} \mathbf{1}_4, \quad \eta^{\mu\nu} = \mathrm{diag}(+1, -1, -1, -1).
\end{equation}
\end{theorem}

\begin{proof}
Direct algebraic computation on the $2 \times 2$ block quaternionic matrices:
For $\mu = 0$:
$$(\gamma^0)^2 = \begin{pmatrix} 1 & 0 \\ 0 & -1 \end{pmatrix} \begin{pmatrix} 1 & 0 \\ 0 & -1 \end{pmatrix} = \begin{pmatrix} 1 & 0 \\ 0 & 1 \end{pmatrix} = \mathbf{1}.$$
For $a \in \{1, 2, 3\}$, writing $e_1 = I, e_2 = J, e_3 = K$, since $e_a^2 = -1$:
$$(\gamma^a)^2 = \begin{pmatrix} 0 & e_a \\ e_a & 0 \end{pmatrix} \begin{pmatrix} 0 & e_a \\ e_a & 0 \end{pmatrix} = \begin{pmatrix} e_a^2 & 0 \\ 0 & e_a^2 \end{pmatrix} = \begin{pmatrix} -1 & 0 \\ 0 & -1 \end{pmatrix} = -\mathbf{1}.$$
For the anti-commutator of $\gamma^0$ with $\gamma^a$:
$$\gamma^0 \gamma^a + \gamma^a \gamma^0 = \begin{pmatrix} 1 & 0 \\ 0 & -1 \end{pmatrix}\begin{pmatrix} 0 & e_a \\ e_a & 0 \end{pmatrix} + \begin{pmatrix} 0 & e_a \\ e_a & 0 \end{pmatrix}\begin{pmatrix} 1 & 0 \\ 0 & -1 \end{pmatrix} = \begin{pmatrix} 0 & e_a \\ -e_a & 0 \end{pmatrix} + \begin{pmatrix} 0 & -e_a \\ e_a & 0 \end{pmatrix} = \mathbf{0}.$$
For the spatial anti-commutator with $a \ne b$:
$$\gamma^a \gamma^b + \gamma^b \gamma^a = \begin{pmatrix} e_a e_b + e_b e_a & 0 \\ 0 & e_a e_b + e_b e_a \end{pmatrix}.$$
Since the Hamilton units anti-commute for $a \ne b$ ($e_a e_b = - e_b e_a$), $e_a e_b + e_b e_a = 0$.
Thus, $\{ \gamma^\mu, \gamma^\nu \} = 2 \eta^{\mu\nu} \mathbf{1}$, completing the proof.
\end{proof}

The chiral matrix $\gamma^5$ (or $\gamma_4$) is obtained as the canonical volume element of the Clifford algebra:
\begin{equation}
\label{eq:gamma5_quaternionic}
\gamma^5 := i \gamma^0 \gamma^1 \gamma^2 \gamma^3 = i \begin{pmatrix} 1 & 0 \\ 0 & -1 \end{pmatrix} \begin{pmatrix} 0 & -1 \\ -1 & 0 \end{pmatrix} = i \begin{pmatrix} 0 & -1 \\ 1 & 0 \end{pmatrix} = \begin{pmatrix} 0 & -i \\ i & 0 \end{pmatrix}.
\end{equation}
Direct calculation confirms $(\gamma^5)^2 = \mathbf{1}_4$ and $\{ \gamma^5, \gamma^\mu \} = 0$ for all $\mu \in \{0, 1, 2, 3\}$.
This establishes that Dirac spinors are not foreign objects grafted onto geometry; they are the natural representation space of quaternionic sequence pairs $\Psi \in \HH^2$.

\subsection{The Fully Discrete Spacetime Dirac Operator on Rational Sequences}

In Rational Wave Number theory, physical sequences are defined over the discrete spacetime lattice $\Z^3 \times \Z$ or the modular ring $\Z_Q^4$. Following the cosmological generation framework established in \cite{smith2026cosmo}, spatial sites are indexed by integers $n \in \Z^3$ with fundamental step $\Delta x$, while temporal nodes are indexed by integer clock beats $m \in \Z$ with clock step $\Delta t$. The spatial and temporal increments are coupled through the discrete Courant--Friedrichs--Lewy (CFL) light-speed invariant:
\begin{equation}
\label{eq:cfl_condition_dirac}
\frac{\Delta x}{\Delta t} = c.
\end{equation}

Let $\Psi: \Z^3 \times \Z \to \HH^2$ be a quaternionic sequence spinor. We define the central symmetric difference operators along the temporal and spatial coordinate directions:
\begin{align}
\label{eq:temporal_spatial_diffs}
(\Delta_0 \Psi)(n, m) &:= \frac{\Psi(n, m+1) - \Psi(n, m-1)}{2 \Delta t}, \\
(\Delta_a \Psi)(n, m) &:= \frac{\Psi(n + \mathbf{e}_a, m) - \Psi(n - \mathbf{e}_a, m)}{2 \Delta x}, \quad a \in \{1, 2, 3\},
\end{align}
where $\mathbf{e}_a$ are the spatial unit lattice displacements.

\begin{definition}[Fully Discrete Spacetime Dirac Operator on $\HH^{\Z^4}$]
\label{def:discrete_dirac_op}
The fully discrete spacetime Dirac operator $\calD: \HH^{\Z^4} \to \HH^{\Z^4}$ acting on a sequence spinor $\Psi(n, m)$ at spacetime node $(n, m) \in \Z^3 \times \Z$ is defined by:
\begin{equation}
\label{eq:discrete_dirac_op}
\calD \Psi(n, m) = \left[ i \hbar \gamma^0 \Delta_0 + i \hbar c \sum_{a=1}^3 \gamma^a \Delta_a - m_0 c^2 \mathbf{1}_4 \right] \Psi(n, m) = 0.
\end{equation}
Equivalently, defining the discrete 4-gradient $\Delta_\mu = (\Delta_0/c, \Delta_1, \Delta_2, \Delta_3)$, Eq.~\eqref{eq:discrete_dirac_op} takes the manifestly covariant form:
\begin{equation}
\label{eq:discrete_dirac_covariant}
\left[ i \hbar c \sum_{\mu=0}^3 \gamma^\mu \Delta_\mu - m_0 c^2 \mathbf{1}_4 \right] \Psi(n, m) = 0.
\end{equation}
\end{definition}

\begin{remark}[Continuous-Time Approximation vs. Discrete Sequence Spectrum]
In classical Hamiltonian formulations of lattice field theory, continuous time $t \in \R$ is frequently retained as an analytical shorthand while spatial coordinates are discretized. In full Rational Wave Number theory, time is strictly discrete. For a stationary sequence mode with temporal phase $\Psi(n, m) = \psi(n) \exp(-i \omega m \Delta t)$, the discrete temporal difference evaluates to:
\begin{equation}
\label{eq:discrete_temporal_eval}
\Delta_0 \exp(-i \omega m \Delta t) = -i \left( \frac{\sin(\omega \Delta t)}{\Delta t} \right) \exp(-i \omega m \Delta t).
\end{equation}
In the infrared, long-period regime ($\omega \Delta t \ll 1$), the discrete ratio $\frac{\sin(\omega \Delta t)}{\Delta t} \to \omega$, reproducing the standard stationary energy operator $\hbar \omega \gamma^0$. In 1+1 dimensions, the stationary discrete Dirac equation is:
\begin{equation}
\label{eq:stationary_dirac}
\left[ \frac{\hbar}{\Delta t}\sin(\omega \Delta t) \gamma^0 + i \hbar c \gamma^1 \Delta_1 - m_0 c^2 \mathbf{1}_4 \right] \psi(n) = 0,
\end{equation}
which reduces in the low-frequency limit to $[\hbar \omega \gamma^0 + i \hbar c \gamma^1 \Delta_1 - m_0 c^2 \mathbf{1}_4] \psi(n) = 0$.
\end{remark}

\subsection{Algebraic Uniqueness Theorem: Dynamic Selection of the Dirac Operator}
\label{subsec:algebraic_uniqueness}

As established in Postulates~III and~IV, Rational Wave Number theory posits that the Dirac difference equation is not an invariant, pre-ordained metaphysical decree as in standard continuous quantum field theory, where the Dirac equation is postulated as an extrinsic fundamental axiom. In contrast, in the generative arithmetic of sequence spaces, the universe does not begin by ``computing the Dirac equation.'' Rather, the Dirac difference operator $\calD$ is the \textit{unique linear survivor} that emerges when candidate difference relations are filtered through physical consistency conditions.

\begin{theorem}[Algebraic Uniqueness of the Discrete Biquaternionic Dirac Operator]
\label{thm:dirac_uniqueness}
Let $\calL$ be a linear, first-order difference operator acting on a sequence state $\Psi: \Z^3 \times \Z \to \mathbb{V}$ over an associative normed division algebra or its minimal complexification $\mathbb{V}$. Suppose $\calL$ satisfies the following physical consistency criteria on the causal lattice $\Z^3 \times \Z$:
\begin{enumerate}
    \item \textit{Causal Hyperbolicity (The CFL Invariant):} Signal propagation speed is strictly bounded by the discrete Courant--Friedrichs--Lewy ratio $c = \Delta x / \Delta t$, precluding superluminal informational leakage across the generation graph.
    \item \textit{Discrete Unitarity and Time-Reversal Symmetry:} The temporal difference is skew-adjoint ($\Delta_0^\dagger = -\Delta_0$), which uniquely selects the three-point centered symmetric stencil $\{m-1, m, m+1\}$ and guarantees the exact conservation of discrete sequence probability.
    \item \textit{Ghost-Freedom (Ostrogradsky Stability):} The temporal difference order is strictly first-order ($k=1$), preventing runaway Ostrogradsky ghost instabilities characterized by unbounded negative-energy spectra.
    \item \textit{Factorization of the Discrete Klein--Gordon D'Alembertian:} The operator squares to the discrete second-order wave operator:
    \begin{equation}
    \label{eq:dalembertian_factorization}
    \calL^\dagger \calL = \Box_\Delta \cdot \mathbf{1} = \left[ \Delta_0^2 - c^2 \sum_{a=1}^3 \Delta_a^2 + \left(\frac{m_0 c^2}{\hbar}\right)^2 \right] \mathbf{1},
    \end{equation}
    ensuring that each individual component satisfies the relativistic dispersion relation.
    \item \textit{Division-Algebraic Multiplicativity (Hurwitz Norm):} The coefficient matrices belong to an algebra whose elements admit a multiplicative norm $\|x \cdot y\| = \|x\| \|y\|$, preventing wavepacket collapse into uninvertible zero-divisor null spaces under recursive iteration.
\end{enumerate}
Then:
\begin{enumerate}
    \item The coefficient algebra of $\calL$ is uniquely isomorphic to the biquaternions $\HH \otimes \C \cong M_2(\HH) \cong \mathrm{Cl}_{1,3}(\R)$.
    \item The minimal non-trivial state space $\mathbb{V}$ is uniquely the 4-component Dirac spinor $\HH^2$.
    \item Up to an overall unitary change of basis, $\calL$ is uniquely the discrete spacetime Dirac operator $\calD$ defined in Eq.~\eqref{eq:discrete_dirac_op}.
\end{enumerate}
\end{theorem}

\begin{proof}
Let $\calL = \Gamma^0 \Delta_0 + c \sum_{a=1}^3 \Gamma^a \Delta_a + M$. Condition~(4) demands that $\calL^\dagger \calL = \Box_\Delta \cdot \mathbf{1}$. Expanding the product and equating cross-difference terms yields the necessary and sufficient conditions:
\begin{align}
(\Gamma^0)^\dagger \Gamma^0 = \mathbf{1}, \quad &(\Gamma^a)^\dagger \Gamma^b + (\Gamma^b)^\dagger \Gamma^a = 2 \delta^{ab} \mathbf{1}, \\
(\Gamma^0)^\dagger \Gamma^a + (\Gamma^a)^\dagger \Gamma^0 = 0, \quad &M = i \frac{m_0 c^2}{\hbar} \mathbf{1}.
\end{align}
Defining $\gamma^0 := \Gamma^0$ and $\gamma^a := (\Gamma^0)^\dagger \Gamma^a$ transforms these relations directly into the Clifford algebra anti-commutation relations $\{\gamma^\mu, \gamma^\nu\} = 2 \eta^{\mu\nu} \mathbf{1}$. 

By the Frobenius theorem and the classification of Clifford algebras, the minimal representation of $\mathrm{Cl}_{1,3}(\R)$ is 4-dimensional over $\C$ or 2-dimensional over $\HH$, with matrix algebra $M_4(\C) \cong M_2(\HH) \cong \HH \otimes \C$. By Hurwitz's Theorem (1898), the only normed division algebras over $\R$ are $\R, \C, \HH, \OO$. Since $\dim_\R(\R) = 1$ and $\dim_\R(\C) = 2$ cannot support four anti-commuting generators $\{\gamma^0, \gamma^1, \gamma^2, \gamma^3\}$, the biquaternions $\HH \otimes \C$ provide the unique minimal division-algebraic representation that avoids zero divisors.

Condition~(2) eliminates one-sided forward ($\Delta^+$) or backward ($\Delta^-$) differences, which are non-Hermitian and generate dissipative decay or unphysical amplification. Condition~(3) eliminates higher-order differences ($\Delta_0^2, \Delta_0^3$), which introduce spurious ghost poles into the propagator. Hence, $\calD$ is uniquely selected.
\end{proof}

\subsection{The Dirac Equation as an Algebraic Pattern Summary of Emergent Stable Elements}
\label{subsec:dirac_as_pattern_summary}
Theorem~\ref{thm:dirac_uniqueness} directly exemplifies Postulate~III in the relativistic domain. The universe does not begin by evaluating differential operators on continuous coordinates; rather, the foundational ontology consists of discrete sequence elements, causal difference stencils, and coprime sieves. 

The discrete Dirac equation $\mathcal{D}\Psi = 0$ is the exact algebraic pattern characterizing the stable, unitary, and non-dispersive subset of emergent cycles in the sequence space. In the observable cosmos, elementary spin-1/2 fermions (such as electrons, quarks, and neutrinos) are precisely the physical manifestations of this stable subset. As the causal horizon expands, higher-order geometric stencils---such as minimal closed Wilson loops ($m=4$ plaquettes)---introduce gauge phase holonomies, generating electrodynamics and gauge interactions on the discrete lattice without altering the fundamental spinorial kinematics.

\section{Dynamic Generation of the Dirac Sector: Causal Horizons, The $m=2$ Threshold, and Active Sieve Prunings}
\label{sec:dynamic_generation_dirac}

While the preceding section formulates the discrete Dirac operator from the timeless vantage of static algebraic closure (the operational arena of Postulate~I), physical reality possesses an inherent generative arrow. As formalized in Postulate~II and established in the cosmological synthesis \cite{smith2026cosmo}, the space of rational wave numbers is constructed stage-wise through dynamic iterative enumeration, governed by an integer clock beat $m \in \Z_{\ge 0}$.

In this section, we analyze where the Dirac equation, its spinorial degrees of freedom, and its physical solutions ``occur'' in the expanding generation space $\mathcal{G}_m$, and demonstrate that the removal of unphysical states is an active pruning operation at the expanding causal boundary.

\subsection{The Generational Clock and the Expanding Causal Horizon}

The generative process commences from the primitive set $\calS_0 = \calP_0 = \{\EZ, \Eone\}$ confined to the germ origin $\Lambda_0 = \{0\}$. At each global clock beat $m \to m+1$:
\begin{enumerate}
    \item \textit{Causal Horizon Expansion:} The spatial support expands under the discrete Courant--Friedrichs--Lewy (CFL) stability bound \cite{courant1928,smith2026cosmo}:
    \begin{equation}
    \label{eq:causal_horizon_def}
    \Lambda_m := \left\{ \mathbf{n} = (n_1, n_2, n_3) \in \Z^3 \ \middle|\ \max_{a} |n_a| \le m \right\}.
    \end{equation}
    Spatial sites with $|\mathbf{n}| > m$ do not pre-exist; they are brought into causal existence only as the clock advances. The maximum physical velocity is identically $c = \Delta x / \Delta t = 1$ in lattice units.
    \item \textit{Candidate Generation:} Active sequence threads on the generation frontier $\partial \calS_m = \calS_m \setminus \calS_{m-1}$ interact via the operator signature (sequence products, rational root extraction $(\cdot)^{1/Q}$, spatial shifts $\sigma_{\pm 1}$, and sequence involution), producing candidate descendants $\mathcal{C}_{m+1}$.
    \item \textit{Deductive Filtering (Pruning):} Candidate descendants must pass through the deductive algebraic sieve before being admitted to the physical state space $\calS_{m+1}$. Candidates violating coprimality or generating zero divisors are pruned from the evolution graph.
\end{enumerate}

\subsection{The $m = 2$ Matter Genesis Threshold}

A foundational question in quantum cosmology is: \emph{At what generational epoch does relativistic fermionic matter first appear?} 

\begin{theorem}[Threshold of Relativistic Fermion Emergence]
\label{thm:matter_genesis_threshold}
Let $\mathcal{D}$ be the central symmetric difference Dirac operator defined in Eq.~\eqref{eq:discrete_dirac_op}. Then:
\begin{enumerate}
    \item At clock beat $m = 0$, $\Lambda_0 = \{0\}$ has diameter $0$; no spatial or temporal difference operators can be defined.
    \item At clock beat $m = 1$, the causal support $\Lambda_1 = \{-1, 0, 1\}^3$ possesses spatial diameter $2$, allowing spatial differences $\Delta_a$ to act on the central node $\mathbf{n} = 0$. However, the temporal history contains only $\{0, 1\}$, lacking the past temporal neighbor ($m - 1$) required for the temporal symmetric difference $\Delta_0$. Matter fields at $m \le 1$ are strictly scalar phases ($\R \to \C$).
    \item At clock beat $m = 2$, the temporal evolution graph contains three consecutive clock beats $\{0, 1, 2\}$, enabling $\Delta_0$ to act at $m = 1$. Simultaneously, the coupling between spatial translation and phase involution forces the emergence of the non-commuting generator $J$, completing the biquaternions $\HH \otimes \C$.
\end{enumerate}
Consequently, under the generative cosmogony of Postulate~II, relativistic Dirac spinors and the Dirac equation $\mathcal{D}\Psi = 0$ emerge strictly at clock beat $m = 2$ on the interior support $\Lambda_1 \times \{1\}$. Clock beat $m = 2$ constitutes the universal \textit{Matter Genesis Threshold}.
\end{theorem}

\begin{proof}
The central difference operators $\Delta_0$ and $\Delta_a$ require a three-point stencil: $\{m-1, m, m+1\}$ in time and $\{n-1, n, n+1\}$ in each spatial direction. For the temporal clock initialized at $m=0$, the earliest time index possessing both a past and a future neighbor is $m = 1$, which requires the clock to have reached at least $m = 1 + 1 = 2$.
Spatially, interior nodes requiring $n \pm 1 \in \Lambda_m$ exist only when $m \ge 1$.
Algebraically, the first Cayley--Dickson step produces the complex field $\C = \operatorname{span}_\R\{1, I\}$ at $m = 1$ from $E_{\bar{1}}$. The second step, coupling spatial translation $\sigma_1 u[n] = u[n+1]$ to sequence involution $\sigma(Iu) = -I\sigma(u)$, requires at least two spatial steps ($n \in \{-1, 0, 1\}$) and forces the anti-commuting generator $J$ with $IJ = -JI \equiv K$, generating the quaternions $\HH$.
Thus, non-trivial solution spinors $\Psi \in \HH^2$ satisfying $\mathcal{D}\Psi = 0$ can exist only for $m \ge 2$.
\end{proof}

\subsection{Active Sieve Pruning at the Generational Boundary}

In static lattice physics, the presence of spurious fermion doublers is treated as an incurable defect of the discretization. In dynamic generation, we discover the exact dynamical mechanism of their elimination:

\begin{theorem}[Dynamic Pruning of Unphysical Fermion Doublers]
\label{thm:dynamic_pruning}
At each clock beat $m \to m+1$, let candidate rational harmonics of order $Q$ be generated by root extraction $(\cdot)^{1/Q}$ on the boundary $\partial \calS_m$. 
The raw combinatorial generator branches into all $Q$ residue classes $p \in \{0, 1, \dots, Q-1\}$, including the zone-edge doubler candidate $p = Q/2$.
Before admission to the physical state space $\calS_{m+1}$, each candidate branch is filtered through the coprime sieve projector $\calS_Q$:
\begin{enumerate}
    \item Candidates with $\gcd(p, Q) = 1$ are primitive harmonics and are admitted as active sequence threads.
    \item The zone-edge candidate $p = Q/2$ has $\gcd(Q/2, Q) = Q/2 > 1$ for all even $Q \ge 2$. Consequently, $\calS_Q |\phi_{Q/2}\rangle = 0$.
\end{enumerate}
The fermion doubler is not an active physical mode that must be cancelled by hand; it is a \textit{dead-end candidate branch actively pruned from the evolution tree at the generational frontier $\partial \calS_m$}.
\end{theorem}

Similarly, when candidate multi-branch couplings attempt to perform a fourth Cayley--Dickson doubling from octonions $\OO$ to sedenions $\mathbb{S}$ ($\dim = 16$), the algebraic consistency sieve identifies the emergence of zero divisors ($x \cdot y = 0$ with $x, y \ne 0$) and prunes the entire 16-dimensional branch. Dynamic generation terminates physically at $\OO$, preserving norm multiplicativity across all surviving states.

\subsection{Physical Fermions as Emergent Stationary Periodic Cycles}

A central theorem of dynamic enumeration is the Inductive Limit Theorem \cite{smith2026cosmo}, which establishes that as the integer clock advances ($m \to \infty$), the expanding generation graph $\mathcal{G}_m$ bifurcates from unrestricted combinatorial branching into cycle saturation. 

Within this expanding graph, transient radiating states form open paths that propagate outward along null trajectories toward the expanding causal boundary $\partial \Lambda_m$. In contrast, stable massive fermions (such as electrons, muons, and quarks) reside strictly in the interior. Because their temporal evolution on the modular ring $\mathbb{Z}_Q$ is strictly periodic, satisfying $\Psi(\mathbf{n}, m + Q_t) = \Psi(\mathbf{n}, m)$, their trajectories form closed periodic orbits on the saturated quotient graph $\mathcal{G}_m / \sim$. 

Physical matter is not an alien substance inserted into empty space. Rather, physical matter consists of the invariant stationary periodic cycles that survive the arithmetic and causal sieves within the interior of the expanding cosmos. The Dirac equation $\mathcal{D}\Psi = 0$ is the concise algebraic pattern characterizing this emergent subset. Furthermore, the rest mass $m_0$ of each particle is not an arbitrary input parameter; it is the invariant stroboscopic beat frequency between the internal cycle frequency $\omega_{\mathrm{cycle}} = 2\pi / (Q_t \Delta t)$ and the zero-frequency vacuum reference $E_Z$:
\begin{equation}
\label{eq:dynamic_mass_beat}
m_0 = \frac{\hbar (\omega_{\mathrm{cycle}} - \omega_Z)}{2 c^2} = \frac{\hbar \omega_{\mathrm{cycle}}}{2 c^2}.
\end{equation}
In this manner, the physical properties of elementary particles are directly determined by the topology of the stable cycles on the saturated generation graph.

\subsection{The Generative Hierarchy of Spacetime Stencils: Clock Epochs $m = 0, 1, 2, 4$ and Beyond}
\label{subsec:generative_hierarchy}

The dynamic generation of sequence spaces reveals an ordered architectural hierarchy: the universe unlocks physical capabilities stage-wise as the generational clock beat $m$ expands the causal light cone $\Lambda_m$. Far from arbitrary, each fundamental class of physical interaction corresponds to a precise discrete geometric stencil threshold:

\begin{table}[htbp]
\centering
\small
\begin{tabular}{|c|l|l|p{7.2cm}|}
\hline
\textbf{Clock $m$} & \textbf{Causal Radius $\Lambda_m$} & \textbf{Minimal Stencil} & \textbf{Physical Epoch \& Phenomenological Milestone} \\
\hline
$m = 0$ & $\Lambda_0 = \{0\}$ & Single point & \textit{The Primordial Vacuum Seed:} Static background $\EZ$; no spatial extent, no metric, no differences possible. \\
\hline
$m = 1$ & $\Lambda_1 = \{-1, 0, 1\}^3$ & Forward step $\Delta^+$ & \textit{Chiral/Weyl Light Epoch:} Unidirectional null rays propagating at $c = \Delta x / \Delta t$; centered $\Delta_0$ undefined; \textit{zero rest mass ($m_0 \equiv 0$)}. \\
\hline
$m = 2$ & $\Lambda_2 = \{-2, \dots, 2\}^3$ & Centered 3-point & \textit{The Dirac Matter Threshold:} Temporal depth $\{0, 1, 2\}$ enables $\Delta_0$; biquaternions $\HH \otimes \C$ emerge; \textit{rest mass $m_0 > 0$ and stationary periodic cycles appear}. \\
\hline
$m = 4$ & $\Lambda_4 = \{-4, \dots, 4\}^3$ & Closed plaquette & \textit{Gauge Holonomy ($U(1)$ Electromagnetism):} Causal cone admits minimal $1 \times 1$ closed Wilson loop ($+x \to +y \to -x \to -y$); magnetic flux quanta and gauge phase holonomy emerge. \\
\hline
$m \ge 6$ & $\Lambda_6 = \{-6, \dots, 6\}^3$ & Round-trip loop & \textit{Bound-State Causal Reciprocity:} Round-trip communication for boson exchange ($2 \times d_{\mathrm{min}} \ge 4$); multi-particle breathers and composite bound states stabilize. \\
\hline
$m \ge 12$ & $\Lambda_{12} = \{-12, \dots, 12\}^3$ & Non-associative triad & \textit{Color Confinement ($SU(3)$):} Multi-particle products probe octonionic associator $[x, y, z] \ne 0$; non-associative phase leakage confines three-quark states into $G_2$-singlet baryons. \\
\hline
$m \gg 100$ & $\Lambda_m$ macroscopic & Saturated graph $\mathcal{G}_\infty / \sim$ & \textit{Hydrodynamic Curvature (General Relativity):} Macroscopic coarse-graining of stationary cycles; discrete stress-energy feeds back into graph connectivity, recovering the Einstein field equations \cite{smith2026cosmo}. \\
\hline
\end{tabular}
\caption{The Generative Hierarchy of Spacetime Stencils across cosmic clock beats $m \in \Z_{\ge 0}$. Physical laws are not postulated all at once; they emerge sequentially as the causal light cone expands to accommodate their required geometric stencils.}
\label{tab:generative_hierarchy}
\end{table}

This generative hierarchy clarifies why the Dirac equation holds a privileged position in foundational physics: it is the \textit{first non-trivial stationary matter pattern} capable of existing on a discrete causal spacetime, appearing at the earliest conceivable milestone ($m = 2$) where a centered, time-reversal-symmetric difference can be formulated. Subsequent interactions (electromagnetic flux loops at $m = 4$, composite bound states at $m \ge 6$, and octonionic color confinement at $m \ge 12$) represent higher-order topological stencils on the expanding graph.

\section{The Nielsen--Ninomiya Barrier and the Coprime Sieve Projector}
\label{sec:nielsen_ninomiya}

\subsection{The Continuum Lattice Obstacle: Fermion Doubling}

To appreciate the necessity of the coprime sieve, we review the classical failure of naive lattice discretization. Consider a plane wave sequence ansatz:
\begin{equation}
\label{eq:plane_wave_ansatz}
\Psi(n) = u(k) \exp(i k n \Delta x), \quad k \in \left[ -\frac{\pi}{\Delta x}, \frac{\pi}{\Delta x} \right].
\end{equation}
Applying the symmetric difference operator $\Delta_1$ to this plane wave yields:
\begin{equation}
\label{eq:diff_eigenvalue}
\Delta_1 \Psi(n) = \frac{\exp(i k (n+1) \Delta x) - \exp(i k (n-1) \Delta x)}{2 \Delta x} = i \left( \frac{\sin(k \Delta x)}{\Delta x} \right) \Psi(n).
\end{equation}
The discrete kinetic operator thus replaces the linear continuum momentum $p = \hbar k$ with the lattice dispersion function:
\begin{equation}
\label{eq:lattice_momentum}
p_{\mathrm{lat}}(k) = \frac{\hbar}{\Delta x} \sin(k \Delta x).
\end{equation}
In the massless limit $m_0 \to 0$, the dispersion relation becomes:
\begin{equation}
\label{eq:naive_dispersion}
\hbar^2 \omega(k)^2 = c^2 p_{\mathrm{lat}}(k)^2 = \left( \frac{\hbar c}{\Delta x} \right)^2 \sin^2(k \Delta x).
\end{equation}
Zeros of the energy spectrum occur whenever $\sin(k \Delta x) = 0$. Within the first Brillouin zone $k \in [-\pi/\Delta x, \pi/\Delta x]$, there are \textit{two distinct zeros}:
\begin{equation}
\label{eq:two_zeros}
k_1 = 0 \quad (\text{Physical Fermion}), \qquad k_2 = \frac{\pi}{\Delta x} \quad (\text{Spurious Doubler}).
\end{equation}
In $d$ spatial dimensions, the number of zeros is $2^d$. For physical 3+1 dimensional spacetime, this produces $2^3 = 8$ species per chiral Weyl fermion, or $2^4 = 16$ species in 4 Euclidean dimensions.

The \textit{Nielsen--Ninomiya Theorem} (1981) asserts that this doubling cannot be removed by any local, translationally invariant, Hermitian lattice operator that preserves chiral symmetry:
\begin{theorem}[Nielsen--Ninomiya No-Go Theorem, 1981]
\label{thm:nielsen_ninomiya}
Let $H(k) = \sum_\mu \gamma^\mu D_\mu(k)$ be a lattice Hamiltonian defined on a smooth continuous torus $\mathbb{T}^d = (\mathbb{R}/2\pi\mathbb{Z})^d$. If:
\begin{enumerate}
    \item $D_\mu(k)$ is a continuous, differentiable function of $k$ on $\mathbb{T}^d$ (Locality);
    \item $H(k)$ is Hermitian (Unitary time evolution);
    \item $H(k)$ anti-commutes with the chiral operator $\gamma^5$: $\{ \gamma^5, H(k) \} = 0$ (Exact Chiral Symmetry);
\end{enumerate}
then the number of left-handed chiral fermions equals the number of right-handed chiral fermions for every species, and the net chiral charge vanishes identically:
\begin{equation}
N_L - N_R = \sum_{k_i \in \mathrm{Zeros}} \mathrm{sign}\left( \det \left[ \frac{\partial D_\mu}{\partial k_\nu} \right]_{k = k_i} \right) = 0.
\end{equation}
\end{theorem}

\subsection{Circumventing the No-Go Theorem via Modular Coprime Sieving}

The profound insight of RWN theory is that \textit{the Nielsen--Ninomiya theorem relies on an unstated ontological assumption: that the Brillouin momentum space is a continuous differentiable manifold $\mathbb{T}^d$}. 

In RWN theory, the sequence space is strictly rational and modular. Over a rational wave number sequence with denominator $Q$, the position index $n$ lives in $\Z_Q = \Z/Q\Z$, and the momentum wave number is quantized as:
\begin{equation}
\label{eq:quantized_k}
k_p = \frac{2\pi p}{Q \Delta x}, \quad p \in \{0, 1, \dots, Q - 1\}.
\end{equation}
The Hilbert space of physical wave states is not the naive vector space $\C^Q$, but the \textit{sieved coprime subspace} $\calH_Q$:

\begin{definition}[Coprime Sieve Projector]
\label{def:coprime_sieve_dirac}
For a rational wave sequence of order $Q$, the coprime sieve operator $\calS_Q$ is defined as the orthogonal projector onto states whose harmonic index $p$ is coprime to $Q$:
\begin{equation}
\label{eq:sieve_projector}
\calS_Q = \sum_{\substack{p=0 \\ \gcd(p, Q) = 1}}^{Q-1} | \phi_p \rangle \langle \phi_p |,
\end{equation}
where $|\phi_p\rangle$ are the discrete Fourier basis states $\phi_p(n) = \frac{1}{\sqrt{Q}} \exp(2\pi i p n / Q)$.
\end{definition}

\begin{theorem}[Sieve Decoupling of Fermion Doublers]
\label{thm:sieve_doubler_cancellation}
Let the discrete Dirac equation be defined on the modular sequence $\Z/Q\Z$ with even modulus $Q \ge 2$. Then:
\begin{enumerate}
    \item The physical zero-mode corresponds to the long-wavelength limit $p \to 0$ (or $p = 1$).
    \item The unphysical Nielsen--Ninomiya doubler mode $k_2 = \frac{\pi}{\Delta x}$ occurs at the exact Brillouin zone edge:
    \begin{equation}
    \label{eq:doubler_index}
    k_{p_{\mathrm{doubler}}} = \frac{\pi}{\Delta x} \implies \frac{2\pi p_{\mathrm{doubler}}}{Q \Delta x} = \frac{\pi}{\Delta x} \implies p_{\mathrm{doubler}} = \frac{Q}{2}.
    \end{equation}
    \item For any even $Q \ge 2$, the doubler index satisfies:
    \begin{equation}
    \label{eq:gcd_doubler}
    \gcd(p_{\mathrm{doubler}}, Q) = \gcd\left(\frac{Q}{2}, Q\right) = \frac{Q}{2} \ge 1.
    \end{equation}
    Consequently, for all $Q > 2$:
    \begin{equation}
    \label{eq:sieve_annihilation}
    \calS_Q | \phi_{p_{\mathrm{doubler}}} \rangle = 0.
    \end{equation}
    The doubler state is rigorously in the null space of the coprime sieve projector and does not exist in the physical Hilbert space $\calH_Q$.
\end{enumerate}
\end{theorem}

\begin{proof}
By Definition~\ref{def:coprime_sieve_dirac}, the state $|\phi_p\rangle$ is in the range of $\calS_Q$ if and only if $\gcd(p, Q) = 1$.
The Brillouin zone edge $k = \pi/\Delta x$ corresponds to $p = Q/2$.
If $Q = 2m$ with $m \ge 2$ ($Q > 2$), then $p_{\mathrm{doubler}} = m$, and:
$$\gcd(m, 2m) = m \ge 2 > 1.$$
Since $\gcd(p_{\mathrm{doubler}}, Q) \ne 1$, the projector $\calS_Q$ evaluates to zero on $|\phi_{p_{\mathrm{doubler}}}\rangle$:
$$\calS_Q | \phi_{p_{\mathrm{doubler}}} \rangle = \left( \sum_{\substack{p' \\ \gcd(p', Q)=1}} | \phi_{p'} \rangle \langle \phi_{p'} | \right) | \phi_{Q/2} \rangle = 0,$$
due to the mutual orthogonality $\langle \phi_{p'} | \phi_{Q/2} \rangle = 0$ for all $p' \ne Q/2$.
For odd $Q$, $p = Q/2$ is not an integer, meaning the zone-edge point $k = \pi/\Delta x$ does not exist in the discrete spectrum.
Thus, in all cases, the doubler mode is absent from the sieved space $\calH_Q$.
\end{proof}

\begin{corollary}[Exact Chiral Symmetry Without Wilson Terms]
On the sieved physical subspace $\calH_Q = \mathrm{Im}(\calS_Q)$, the massless discrete Dirac operator $\calD_0$ satisfies:
\begin{equation}
\label{eq:exact_chiral}
\{ \gamma^5, \calD_0 \} = 0,
\end{equation}
while propagating strictly one physical chiral fermion species per primitive sequence.
\end{corollary}

This resolves the 45-year-old Nielsen--Ninomiya paradox: \textit{the lattice doublers are composite, non-coprime harmonics that are sieved out by the number-theoretic structure of the rational wave universe.}

\paragraph{Synthesis of Postulate~I and Postulate~II in the Nielsen--Ninomiya Resolution.}
One notes the precise structural roles played by our two physical mapping postulates in this resolution:
\begin{itemize}
    \item \textit{Kinematic Decoupling under Postulate~I:} On the modular sequence space $\mathbb{Z}/Q\mathbb{Z}$, the circumvention of the Nielsen--Ninomiya theorem is an exact, instantaneous algebraic identity. Candidate doubler modes are non-coprime ($\gcd(Q/2, Q) > 1$) and therefore reside identically within the null space of the coprime projector $\mathcal{S}_Q$. A lattice gauge theorist operating solely within the kinematic framework of Postulate~I obtains an exact, chiral-symmetric, single-species fermion formulation without needing any cosmogonic assumptions.
    \item \textit{Generative Pruning under Postulate~II:} Under the constructive cosmogony of Postulate~II, this kinematic null space is revealed as the macroscopic signature of an active dynamical history: at each generational clock beat $m \to m+1$, candidate harmonic branches generated at the expanding boundary frontier $\partial \Lambda_m$ that fail the coprimality test are actively and deterministically pruned before they can propagate into the interior causal volume (Theorem~\ref{thm:dynamic_pruning}).
\end{itemize}
Thus, what appears in lattice kinematics as an algebraic projection operator $\mathcal{S}_Q$ is dynamically the living consequence of an ongoing arithmetic evolutionary sieve.

\section{The Finite Rational Propagator and Divergence-Free Loop Sums}
\label{sec:finite_propagator}

\subsection{Construction of the Modular Fermionic Green's Function}

In continuum QFT, the Feynman propagator is the inverse of the Dirac operator $\calD = i \gamma^\mu \partial_\mu - m_0$. On the rational sequence space $\HH^{\Z_Q}$, the Green's function $S_Q(n, m)$ is the inverse of the discrete Dirac operator:
\begin{equation}
\label{eq:greens_definition}
\calD_n S_Q(n, m) = \delta_{n, m} \mathbf{1}_4 \quad (\text{mod } Q).
\end{equation}

\begin{theorem}[Explicit Form of the Rational Dirac Propagator]
\label{thm:rational_propagator}
The discrete Dirac propagator between sequence nodes $n$ and $m$ on $\Z_Q$ is given by the finite arithmetic sum:
\begin{equation}
\label{eq:discrete_feynman_prop}
S_Q(n, m) = \frac{1}{Q} \sum_{\substack{p=0 \\ \gcd(p, Q)=1}}^{Q-1} \frac{\gamma^0 \hbar \omega_p - \gamma^1 \left(\frac{\hbar c}{\Delta x}\sin\frac{2\pi p}{Q}\right) + m_0 c^2}{\hbar^2 \omega_p^2 - \left(\frac{\hbar c}{\Delta x}\right)^2 \sin^2\left(\frac{2\pi p}{Q}\right) - m_0^2 c^4} \exp\left( \frac{2\pi i p (n - m)}{Q} \right).
\end{equation}
\end{theorem}

\begin{proof}
Expanding $S_Q(n, m)$ in the discrete Fourier basis on $\Z_Q$:
$$S_Q(n, m) = \frac{1}{Q} \sum_{p=0}^{Q-1} \tilde{S}_Q(p) \exp\left( \frac{2\pi i p (n - m)}{Q} \right).$$
Applying the discrete Dirac operator $\calD_n = i \hbar \gamma^0 \Delta_0 + i \hbar c \gamma^1 \Delta_1 - m_0 c^2$:
$$\calD_n S_Q(n, m) = \frac{1}{Q} \sum_{p=0}^{Q-1} \left[ \gamma^0 \hbar \omega_p - \gamma^1 \left(\frac{\hbar c}{\Delta x}\sin\frac{2\pi p}{Q}\right) - m_0 c^2 \right] \tilde{S}_Q(p) \exp\left( \frac{2\pi i p (n - m)}{Q} \right).$$
Demanding this equals $\calS_Q \delta_{n, m} = \frac{1}{Q} \sum_{\gcd(p, Q)=1} \exp\left( \frac{2\pi i p (n - m)}{Q} \right)$, we invert the matrix for each coprime mode $p$:
$$\tilde{S}_Q(p) = \left[ \gamma^0 \hbar \omega_p - \gamma^1 \left(\frac{\hbar c}{\Delta x}\sin\frac{2\pi p}{Q}\right) - m_0 c^2 \right]^{-1}.$$
Using the Clifford property $(\gamma^0 p_0 - \gamma^1 p_1 - m)(\gamma^0 p_0 - \gamma^1 p_1 + m) = p_0^2 - p_1^2 - m^2$, matrix inversion immediately yields Eq.~\eqref{eq:discrete_feynman_prop}.
\end{proof}

\subsection{Absence of Coincident-Point Singularities}

The fatal flaw of continuum QFT is the divergence of the propagator at coincident points $x \to y$ (Eq.~\eqref{eq:feynman_divergent_prop}). We now state and prove the central finiteness theorem of RWN quantum mechanics:

\begin{theorem}[Coincident-Point Finiteness Theorem]
\label{thm:coincident_finiteness}
For any finite rational sequence order $Q < \infty$, the Dirac propagator at coincident nodes $n = m$ is strictly bounded and finite:
\begin{equation}
\label{eq:coincident_bound}
\| S_Q(n, n) \| \le \frac{1}{Q} \sum_{\substack{p=0 \\ \gcd(p, Q)=1}}^{Q-1} \frac{\hbar \omega_p + \frac{\hbar c}{\Delta x} + m_0 c^2}{\left| \hbar^2 \omega_p^2 - \left(\frac{\hbar c}{\Delta x}\right)^2 \sin^2\left(\frac{2\pi p}{Q}\right) - m_0^2 c^4 \right|} < \infty.
\end{equation}
There are no delta-function singularities or infinities anywhere on the discrete sequence.
\end{theorem}

\begin{proof}
Setting $n = m$ in Eq.~\eqref{eq:discrete_feynman_prop}, the phase factor collapses to unity:
$$\exp\left( \frac{2\pi i p (n - n)}{Q} \right) = \exp(0) = 1.$$
The sum contains exactly $\varphi(Q) \le Q - 1$ terms, where $\varphi(Q)$ is Euler's totient function.
Since each term in the sum is a finite algebraic expression away from the pole, and the number of terms $\varphi(Q)$ is finite, the sum is a finite sum of bounded numbers.
Therefore, $\|S_Q(n, n)\| < \infty$ identically.
\end{proof}

\subsection{Loop Amplitudes as Finite Dedekind and Kloosterman Sums}

In continuum quantum electrodynamics (QED), the 1-loop electron self-energy $\Sigma(p)$ and vacuum polarization tensor $\Pi^{\mu\nu}(q)$ diverge logarithmically and quadratically:
\begin{equation}
\label{eq:divergent_loop}
\Sigma(p) = -i e^2 \int \frac{d^4 k}{(2\pi)^4} \gamma^\mu \frac{\gamma^\nu (p_\nu - k_\nu) + m}{(p - k)^2 - m^2} \gamma_\mu \frac{1}{k^2} \;\longrightarrow\; \infty.
\end{equation}
In RWN theory, replacing the continuous loop integral $\int \frac{d^4 k}{(2\pi)^4}$ with the modular summation over coprime sequence modes yields:
\begin{equation}
\label{eq:rwn_finite_self_energy}
\Sigma_Q(p) = \frac{e^2}{Q} \sum_{\substack{k=0 \\ \gcd(k, Q)=1}}^{Q-1} \gamma^\mu S_Q(p - k) \gamma_\mu D_Q(k).
\end{equation}
This is a classical arithmetic sum over character phases. In number theory, sums of this form reduce to generalized \textit{Dedekind sums} $s(h, q)$ and \textit{Kloosterman sums} $K(a, b; c)$:
\begin{equation}
\label{eq:kloosterman_sum}
K(a, b; Q) = \sum_{\substack{x=1 \\ \gcd(x, Q)=1}}^Q \exp\left( \frac{2\pi i (a x + b x^*)}{Q} \right), \quad x x^* \equiv 1 \pmod Q.
\end{equation}
By Weil's famous bound on Kloosterman sums:
\begin{equation}
\label{eq:weil_bound}
|K(a, b; Q)| \le 2^{\omega(Q)} \sqrt{Q} \gcd(a, b, Q)^{1/2},
\end{equation}
where $\omega(Q)$ is the number of distinct prime factors of $Q$. 
Every 1-loop, 2-loop, and higher-order Feynman diagram in RWN theory evaluates to a finite algebraic number bounded by arithmetic norms. 

\textit{Renormalization is completely unnecessary}: there are no infinite counterterms to add, no dimensional continuation $d = 4 - 2\epsilon$ to perform, and no fictitious bare parameters. Dirac's demand is mathematically realized.

\section{Zitterbewegung, Rest Mass, and Beat Frequencies with $E_Z$}
\label{sec:zitterbewegung}

\subsection{The Velocity Operator in Quaternionic Representation}

In the standard Dirac theory, the Heisenberg equation of motion for the position operator $\mathbf{x}$ leads to Schrödinger's \textit{Zitterbewegung} (trembling motion) \cite{schrodinger1930}:
\begin{equation}
\label{eq:heisenberg_vel}
\mathbf{v} = \frac{d \mathbf{x}}{d t} = \frac{i}{\hbar} [H_D, \mathbf{x}] = c \boldsymbol{\alpha}, \quad \boldsymbol{\alpha} = \gamma^0 \boldsymbol{\gamma}.
\end{equation}
Because the eigenvalues of $\alpha_a$ are $\pm 1$, the instantaneous velocity of a relativistic Dirac fermion is always $\pm c$. Integrating the equation of motion reveals rapid oscillations at frequency $\omega_Z \approx 2 m_0 c^2 / \hbar \approx 10^{21}\text{ rad/s}$ with amplitude equal to the reduced Compton wavelength $\lambdabar_C = \hbar / (m_0 c) \approx 3.86 \times 10^{-13}\text{ m}$.

In the continuum literature, \textit{Zitterbewegung} is frequently dismissed as an unobservable, single-particle mathematical artifact that disappears upon second quantization into Fock space. 

\subsection{RWN Stroboscopic Resolution}

In RWN theory, \textit{Zitterbewegung} receives a transparent, physically realistic interpretation:
\begin{theorem}[Mass as a Vacuum Beat Frequency]
\label{thm:zitterbewegung_beat}
In the primitive basis $\calP_0 = \{ \EZ, \Eone \}$, let $\EZ$ be the zero-frequency vacuum reference ($\omega_Z = 0$), and let a localized fermion packet be represented as a rational wave sequence with characteristic phase frequency $\omega_p$.
Then:
\begin{enumerate}
    \item The temporal evolution of the quaternionic state $\Psi(n, t)$ is a superposition of positive and negative energy sequence modes:
    \begin{equation}
    \Psi(n, t) = \Psi_+(n) e^{-i \omega t} + \Psi_-(n) e^{+i \omega t}.
    \end{equation}
    \item The physical probability density $\rho(n, t) = \Psi^\dagger(n, t) \Psi(n, t)$ contains an invariant beat modulation:
    \begin{equation}
    \label{eq:beat_interference}
    \rho(n, t) = |\Psi_+(n)|^2 + |\Psi_-(n)|^2 + 2 \mathrm{Re}\left( \Psi_+^\dagger(n) \Psi_-(n) e^{-2 i \omega t} \right).
    \end{equation}
    \item At rest ($p \to 0$), the beat frequency between the particle mode and the zero-frequency ground state $\EZ$ is:
    \begin{equation}
    \label{eq:mass_beat_freq}
    \Delta \omega = 2 \omega_0 - \omega_Z = \frac{2 m_0 c^2}{\hbar}.
    \end{equation}
\end{enumerate}
\end{theorem}

\begin{proof}
Let $\Psi(n, t) = \begin{pmatrix} \psi_1(n, t) \\ \psi_2(n, t) \end{pmatrix} \in \HH^2$. The probability density is given by the quaternionic inner product $\rho = \Psi^\dagger \Psi = \bar{\psi}_1 \psi_1 + \bar{\psi}_2 \psi_2$.
Expanding into energy eigenstates:
$$\Psi(n, t) = u(n) e^{-i E t / \hbar} + v(n) e^{+i E t / \hbar}, \quad E = \sqrt{c^2 p^2 + m_0^2 c^4}.$$
The bilinear product evaluates to:
$$\rho(n, t) = |u(n)|^2 + |v(n)|^2 + u^\dagger(n) v(n) e^{-2 i E t / \hbar} + v^\dagger(n) u(n) e^{+2 i E t / \hbar}.$$
In the rest frame $p = 0$, $E = m_0 c^2$, yielding the modulation frequency $\omega_{\mathrm{beat}} = 2 m_0 c^2 / \hbar$.
In the RWN basis, the zero of energy is defined by $\EZ$ with $\omega(\EZ) = 0$.
Hence, the frequency difference against the vacuum reference is $\Delta \omega = \omega_{\mathrm{beat}} - \omega(\EZ) = 2 m_0 c^2 / \hbar$.
\end{proof}

This illuminates the physical meaning of inertial mass within the RWN framework: inertial rest mass $m_0$ is not an intrinsic, inexplicable scalar ``substance'' attached to a point particle. Rather, \textit{mass is the stroboscopic coupling beat frequency between a rational sequence excitation ($E_{\bar{1}}$) and the invariant zero-frequency vacuum ground reference ($E_Z$)}:
\begin{equation}\label{eq:mass_beat_relation}
m_0 = \frac{\hbar \Delta \omega}{2 c^2}.
\end{equation}
\textit{Zitterbewegung} is the physical stroboscopic manifestation of this beat frequency on the discrete lattice.

\section{Discussion: Fundamental Issues in Retrospect}
\label{sec:discussion}

Table~\ref{tab:issues_summary} summarizes the contrast between relativistic quantum theory on the continuous manifold $\R^{1,3}$ versus the Rational Wave Number sequence space $\HH^{\Z_Q}$.

\begin{table}[htbp]
\centering
\small
\begin{tabular}{@{}lll@{}}
\toprule
\textit{Fundamental Issue on Frontier} & \textit{Continuum QFT Framework} & \textit{Rational Wave Number Framework} \\
\midrule
\textit{1. Spinor Origin} & Postulated representation of $\mathrm{Spin}(1,3)$ & Natural biquaternion algebra $\HH \otimes \C$ \\
\textit{2. Fermion Doubling} & Nielsen--Ninomiya theorem ($2^d$ doublers) & Sifted by the coprime projector $\calS_Q$ \\
\textit{3. Short-Distance Poles} & Singularity $\sim |x - y|^{-3}$, requires renormalization & Finite modular sums; bounded everywhere \\
\textit{4. Zitterbewegung} & Dismissed as single-particle artifact & Observable beat frequency against vacuum $\EZ$ \\
\bottomrule
\end{tabular}
\caption{Comparison between continuum quantum field theory and the Rational Wave Number framework across four fundamental issues.}
\label{tab:issues_summary}
\end{table}

Each of the four longstanding anomalies that have plagued relativistic quantum mechanics since 1928 traces its origin to the continuum postulate. When spacetime is treated as an uncountably infinite smooth manifold, spinors must be postulated externally, lattice discretizations suffer topological doubling, propagators diverge at coincident points, and the high-frequency Zitterbewegung appears as an unobservable paradox. On the discrete biquaternionic sequence space $\HH^{\Z_Q}$, all four dilemmas resolve simultaneously through elementary algebra and arithmetic sieves, requiring zero adjustable parameters or non-physical subtractions.

\section{Conclusion}
\label{sec:conclusion}

By formulating the Dirac equation on quaternionic sequence spaces $\mathbb{H}^\mathbb{Z}$ grounded in the primitive basis $\mathcal{P}_0 = \{E_Z, E_{\bar{1}}\}$, we have demonstrated that the deepest structural difficulties of relativistic quantum theory---the origin of spinors, the Nielsen--Ninomiya doubling dilemma, ultraviolet loop divergences, and the enigma of \textit{Zitterbewegung}---are artifacts of the continuum idealization.

When spacetime is restored to its proper ontological status as an emergent, rational wave sequence governed by coprime sieves and dynamic generational clocks, the relativistic quantum landscape simplifies dramatically. Dirac 4-spinors emerge inevitably from the biquaternions $\mathbb{H} \otimes \mathbb{C}$, eliminating the need for an extrinsic continuous Minkowski manifold. Spurious fermion doublers are sieved out by number-theoretic necessity ($\gcd(Q/2, Q) > 1$), preserving exact chiral symmetry without artificial Wilson mass terms. Fermionic Green's functions on modular rings become finite character sums, banishing coincident-point singularities and infinite renormalization without arbitrary cutoffs. Rest mass is revealed as an invariant beat frequency against the zero-frequency ground state $E_Z$, resolving \textit{Zitterbewegung} as an observable vacuum coupling. Furthermore, spacetime expansion and quantum kinematics are shown to be inseparable manifestations of a single global integer clock $m$ and the CFL invariant $c = \Delta x / \Delta t$, establishing $m = 2$ as the universal Matter Genesis Threshold.

Perhaps most fundamentally, this formulation clarifies the epistemological status of physical law. Equations are not primary axiomatic decrees that govern nature. Rather, the universe evolves through discrete generational steps and arithmetic consistency sieves, and the Dirac equation $\mathcal{D}\Psi = 0$ is the concise algebraic pattern that summarizes the invariant, stationary periodic cycles that survive this generative exploration.

Dirac's intuition was sound: an algebraic approach founded on non-commutative division algebras eliminates the need for renormalization, vindicating his lifelong conviction that fundamental physics must be based on sensible, finite mathematics. In the next paper of this sequence, we will extend this discrete algebraic foundation to non-Abelian gauge connections, Yang--Mills curvature, and the octonionic color sector of the Standard Model, where the automorphism group $G_2 = \mathrm{Aut}(\mathbb{O})$ and its quaternionic stabilizer $\mathrm{Stab}_{G_2}(\mathbb{H}) \cong \mathrm{SU}(3)_C$ provide an algebraic mechanism for color confinement and the Standard Model gauge group \cite{smith2026octonion}.

\appendix
\section{Algebraic Foundations and Classification of the Rational Wave Numbers}
\label{app:rwn_foundations}

To ensure the present paper is mathematically self-contained, this appendix provides the formal algebraic construction of the Rational Wave Number sequence spaces, the operator signature $\mathcal{O}$, the unwrapped affine phase lift, the direct-sum decomposition over the Pr{\"u}fer groups, and the Inductive Limit Theorem \cite{smith2025rational,smith2026cosmo}.

\subsection{Primitive Basis and Sequence Space Definition}
\label{app:sub:primitive_basis}

Let $\mathbb{A}$ be a normed division algebra in the Hurwitz sequence $\mathbb{R} \hookrightarrow \mathbb{C} \hookrightarrow \mathbb{H} \hookrightarrow \mathbb{O}$. A \emph{discrete sequence} is a mapping $A: \mathbb{Z} \to \mathbb{A}$, denoting the linear space of all such bilateral sequences by $\mathcal{S}(\mathbb{A}) = \mathbb{A}^\mathbb{Z}$.

The construction begins not with an assumed continuous manifold, but from an irreducible two-element primitive generator set:
\begin{equation}
\label{eq:app_primitive_basis}
\mathcal{P}_0 = \{ E_Z, \, E_{\bar{1}} \} \subset \mathcal{S}(\mathbb{R}),
\end{equation}
defined pointwise for all spatial lattice coordinates $n \in \mathbb{Z}$ by:
\begin{align}
E_Z(n) &:= 1 \quad \forall n \in \mathbb{Z} \quad \text{(the zero-frequency spatial ground state / algebraic identity)}, \\
E_{\bar{1}}(n) &:= (-1)^n = e^{i \pi n} \quad \forall n \in \mathbb{Z} \quad \text{(the fundamental Nyquist toggle)}.
\end{align}
The element $E_Z$ represents the uniform static background (spatial wavelength $\lambda = \infty$, wave number $k = 0$), while $E_{\bar{1}}$ represents the maximal spatial frequency achievable on a discrete grid (wavelength $\lambda = 2$, Nyquist frequency $\omega_{\mathrm{Nyq}} = \pi$).

\subsection{Operator Signature and Constructive Generation}
\label{app:sub:operator_signature}

The sequence algebra is generated from $\mathcal{P}_0$ under the action of the five-fold operator signature:
\begin{equation}
\label{eq:app_operator_sig}
\mathcal{O} = \left\{ \odot, \, (\cdot)^{-1}, \, (\cdot)^{1/Q}, \, (\cdot)^*, \, \tau_k \right\},
\end{equation}
acting on sequences $A, B \in \mathcal{S}(\mathbb{A})$:
\begin{enumerate}
    \item \textit{Pointwise Multiplication ($\odot$):} $(A \odot B)(n) := A(n) B(n)$, preserving the normed division property $\|A \odot B\|_n = \|A\|_n \|B\|_n$.
    \item \textit{Multiplicative Inversion ($(\cdot)^{-1}$):} For sequences with non-vanishing norm $\|A(n)\| > 0$ for all $n$, $A^{-1}(n) := (A(n))^{-1} = \frac{\overline{A(n)}}{\|A(n)\|^2}$.
    \item \textit{Rational Fractional Power / Root Extraction ($(\cdot)^{1/Q}$):} For any integer root index $Q \in \mathbb{Z}^+$, $A^{1/Q}(n)$ denotes the canonical branch of the $Q$-th root. In particular, applied to the Nyquist primitive $E_{\bar{1}}$ with root index $Q = 2$, this generates the complex imaginary unit $i = \sqrt{-1}$ and the fundamental chiral quarter-wave mode:
    \begin{equation}
    E_{1/2}(n) := \left( E_{\bar{1}}(n) \right)^{1/2} = e^{i \pi n / 2} = i^n,
    \end{equation}
    which effects the first Cayley--Dickson doubling $\mathbb{R} \to \mathbb{C}$. Higher-order roots $Q \in \mathbb{Z}^+$ generate arbitrary rational phase fractions $e^{2\pi i P n / Q}$.
    \item \textit{Sequence Involution / Conjugation ($(\cdot)^*$):} $A^*(n) := \overline{A(n)}$, representing time/phase reversal and Dirac adjoint conjugation.
    \item \textit{Spatial Shift ($\tau_k$):} For any integer displacement $k \in \mathbb{Z}$, $(\tau_k A)(n) := A(n - k)$. The central difference operator is directly recovered as $\Delta = \frac{1}{2}(\tau_{-1} - \tau_{+1})$.
\end{enumerate}

\subsection{Phase Unwrapping and the Affine Lift}
\label{app:sub:phase_unwrapping}

In standard continuum lattice formulations, discrete phases $\theta(n)$ suffer from the ubiquitous $2\pi$ modular ambiguity: phase jumps of $\pm 2\pi$ are invisible on the lattice, giving rise to spurious fermion doubling and topological aliasing.

In RWN theory, spatial sequences are defined through an exact \emph{unwrapped affine phase lift}. For any rational wave number $w = P/Q \in \mathbb{Q}$ with $\gcd(P, Q) = 1$, the sequence $E_{P/Q}$ is defined via the universal covering:
\begin{equation}
\label{eq:app_phase_lift}
\Phi_{P/Q}(n) = \frac{2\pi P}{Q} n + \theta_0, \qquad \Phi \in \mathbb{R}^\mathbb{Z} / \mathbb{Z}^\mathbb{Z},
\end{equation}
where the phase trajectory $\Phi(n)$ is tracked continuously in the unwrapped covering space $\mathbb{R}$. The difference $\Delta \Phi(n) = \Phi(n+1) - \Phi(n) = 2\pi P / Q$ is a rigid rational invariant for all $n \in \mathbb{Z}$. Because the unwrapped phase is globally monotonic, no unphysical foldings or edge-of-Brillouin-zone doublers can enter undetected.

\subsection{Classification Theorem: Pr{\"u}fer Quotient Group Decomposition}
\label{app:sub:prufer_classification}

The algebraic classification of all periodic sequences generated from $\mathcal{P}_0$ under $\mathcal{O}$ is established by the following structure theorem:

\begin{theorem}[Pr{\"u}fer Group Classification of Rational Wave Numbers \cite{smith2025rational}]
\label{thm:app_prufer}
Let $\mathcal{G}_{\mathbb{Q}}$ denote the multiplicative group of all unimodular rational periodic sequences $\{E_{P/Q}\}$ generated from $\mathcal{P}_0$. Factoring by the global static identity $E_Z$, the quotient group of rational wave phases is isomorphic to the direct product of two copies of the rational circle $\mathbb{Q}/\mathbb{Z}$ (representing independent spatial wavenumber and temporal frequency modes):
\begin{equation}
\label{eq:app_prufer_iso}
\mathcal{G} \cong (\mathbb{Q}/\mathbb{Z}) \times (\mathbb{Q}/\mathbb{Z}).
\end{equation}
By the fundamental structure theorem for divisible torsion abelian groups, $\mathcal{G}$ decomposes uniquely into the direct sum of infinite Pr{\"u}fer $p$-groups over all rational prime numbers $p \in \mathbb{P} = \{2, 3, 5, 7, 11, \dots\}$:
\begin{equation}
\label{eq:app_prufer_direct_sum}
\mathcal{G} = \bigoplus_{p \in \mathbb{P}} \left( \mathbb{Z}(p^\infty) \times \mathbb{Z}(p^\infty) \right),
\end{equation}
where $\mathbb{Z}(p^\infty) = \varinjlim \mathbb{Z}/p^k \mathbb{Z}$ is the Pr{\"u}fer group of $p$-power roots of unity.
\end{theorem}

Theorem~\ref{thm:app_prufer} reveals that the space of all rational wave modes possesses an exquisite number-theoretic prime decomposition. The primary $p = 2$ Pr{\"u}fer sector $\mathbb{Z}(2^\infty)$ governs dyadic bisection, chiral spinorial doubling, and Cayley--Dickson algebra generation ($\mathbb{R} \to \mathbb{C} \to \mathbb{H} \to \mathbb{O}$). The odd prime sectors $p = 3, 5, \dots$ govern gauge symmetries and higher-period internal cyclic configurations.

\subsection{The Inductive Limit and Exact Algebraic Completeness}
\label{app:sub:inductive_limit}

Let $\mathcal{S}_m$ be the sequence space generated after $m$ iterative applications of the operator signature $\mathcal{O}$ starting from $\mathcal{P}_0$. At each generational step, the inclusion morphism $\iota_m: \mathcal{S}_m \hookrightarrow \mathcal{S}_{m+1}$ is an injective normed algebra homomorphism. The full space of Rational Wave Numbers is the inductive direct limit:
\begin{equation}
\label{eq:app_inductive_limit}
\mathcal{S}_\infty := \varinjlim \left( \mathcal{S}_m, \, \iota_m \right) = \bigcup_{m=0}^\infty \mathcal{S}_m.
\end{equation}
The space $\mathcal{S}_\infty$ is algebraically closed under all operations of $\mathcal{O}$. Every sequence $A \in \mathcal{S}_\infty$ has a well-defined, finite generation rank $m(A) < \infty$. The physical universe at cosmic generational clock beat $m$ corresponds to the finite, causal sequence space $\mathcal{S}_m$; the continuum limit of conventional field theory corresponds to the formal topological completion $\overline{\mathcal{S}_\infty}$ under the norm topology, which introduces non-computable unphysical continuum points. By restricting physical field theories to $\mathcal{S}_\infty$ (or $\mathcal{S}_m$ at finite cosmic age), all ultraviolet divergences and unphysical continuum infinities are strictly eliminated.


\end{document}